%% file: ex_article.tex
\documentclass[review,onefignum,onetabnum]{siamart171218}

\input{ex_shared}

\ifpdf
\hypersetup{
  pdftitle={Ginzburg-Landau Spiral Waves in Circular and Spherical Geometries},
  pdfauthor={Jia-Yuan Dai}
}
\fi

\begin{document}

\nolinenumbers

\maketitle

\begin{abstract}
  We prove the existence of $m$-armed spiral wave solutions for the complex Ginzburg-Landau equation in the circular and spherical geometries. We establish a new global bifurcation approach and generalize the results of existence for rigidly-rotating spiral waves. Moreover, we prove the existence of two new patterns: frozen spirals in the circular and spherical geometries, and 2-tip spirals in the spherical geometry.
\end{abstract}

\begin{keywords}
 complex Ginzburg-Landau equation, $m$-armed spiral waves, global bifurcation.
\end{keywords}

\begin{AMS}
  35Q56, 37G40
\end{AMS}

\section{Introduction}

The aim of this paper is to establish a global bifurcation approach and bring new existence results in the study of Ginzburg-Landau vortex and spiral wave solutions. 

\subsection{Literature review}

Pattern formation of spiral waves on two-dimensional manifolds has been observed in physiological, chemical, and biological models. In physiology, spiral waves of electricity appear on heart tissues during cardiac arrhythmia and life-threatening fibrillation; see \cite{Daal73, WiRo46}. In chemistry, the diffusive Belousov-Zhabotinsky reaction triggers intricate spiral waves on a petri dish; see \cite{Beetal97, Wi72}. In biology, spiral waves arise during aggregation of slime mold via chemotactic movement; see \cite{FaLe98}.

Spiral waves appear both in excitable and oscillatory media; see \cite{FiSc03, Mu03}. There are two alternative methods to model the dynamics of spiral waves. The first method uses reaction-diffusion systems, in which the reaction kinetics produce the effect of media. Spiral waves are mostly triggered by a mechanism composed of the Hopf bifurcation and symmetry breaking; see \cite{Sc98, SiMa11}. The second method aims at a geometric description of spiral wave fronts by mean curvature flows of curves, and the effect of media is implemented in the normal velocity of the curve; see the survey \cite{Mi03}. Note that the second method is sometimes called the kinematic theory and especially applicable for excitable media such as the diffusive Belousov-Zhabotinsky reaction and FitzHugh-Nagumo model; see \cite{Fietal06, Ke86, MiZy91}.

Although the mechanism that triggers spiral waves is understood, only a few rigorous results of existence are available. One of the popular models that yields spiral waves is the cubic supercritical Ginzburg-Landau equation on $\mathbb{R}^2$,
\begin{equation} \label{gle-1}
\partial_t \Psi  = \Delta_{\mathbb{R}^2} \Psi + ( 1 - |\Psi|^2 - i\, \beta \, |\Psi|^2 ) \, \Psi.
\end{equation}
Here the unknown $\Psi$ is complex valued, and $\beta \in \mathbb{R}$ is a prescribed kinetic parameter.

The significance of (\ref{gle-1}) is twofold. First, in science it plays a central role in the theory of nonlinear hydrodynamics and condensed matter physics; see \cite{Arkr02, Pi06}. Second, in mathematics it is a normal form for general parameter-dependent PDEs near the Hopf instability; see \cite{Mi02, Sc98}. As a normal form, (\ref{gle-1}) faithfully approximates the dynamics of the original underlying PDEs in a qualitative manner. Hence the existence of Ginzburg-Landau spiral waves on $\mathbb{R}^2$ has been served as a theoretical justification for numerical evidence of spiral waves governed by many other models, for instance, by
cyclic competition of species described by the rock-paper-scissors game; see \cite{Reetal08}.

The main feature that (\ref{gle-1}) allows one to prove the existence of spiral waves is the \textit{global $S^1$-equivariance} (or \textit{global gauge symmetry}; see \cite{Arkr02}): For each $\vartheta \in S^1 \cong \mathbb{R}/ 2 \pi \mathbb{Z}$, 
\begin{equation} \label{glo-gau}
\Psi \mbox{    is a solution of (\ref{gle-1}) if and only if   } e^{i \vartheta}\, \Psi \mbox{    is a solution}.
\end{equation}
This symmetry is closely related to the appearance of rotating wave solutions; see \cite{Da17} Section 3.4. More specifically, it allows one to pursue the following \textit{$m$-armed spiral Ansatz} in polar coordinates $(s,\varphi)$ on $\mathbb{R}^2$:
\begin{equation} \label{spians-1}
\Psi(t,s, \varphi) = e^{-i\Omega t}  \left(A(s) \, e^{i p(s)} \right) e^{im\varphi};
\end{equation}
see \cite{Coetal78, Gr80, Ha82, KoHo81}. Notice that most excitable media such as the FitzHugh-Nagumo model do not admit such a symmetry.

Once a nontrivial solution of (\ref{gle-1}) in the form of (\ref{spians-1}) exists, its associated pattern consists of a tip and a spiral-like shape. The tip is a jump discontinuity of the phase field, 
\begin{equation}
P(t,s,\varphi) := -\Omega t + p(s) + m\varphi,
\end{equation} 
and resides at $s = 0$. The shape is seen as the zero contour of the phase field on $\mathbb{R}^2$,
\begin{equation}
\{ \left( s\cos(\varphi(t, s)), s \sin(\varphi(t, s)) \right): P(t, s, \varphi(t,s)) = 0 \},
\end{equation}
and it exhibits a rigidly-rotating spiral by our Ansatz (\ref{spians-1}).

Substituting the spiral Ansatz (\ref{spians-1}) into (\ref{gle-1}) results in the ODEs for the amplitude $A(s)$ and the phase derivative $p'(s)$:
\begin{align} \label{ode-1}
& A'' + \frac{1}{s}  A' - \frac{m^2}{s^2} A -A \, (p')^2 + (1 - A^2) \, A  = 0,
\\ \label{ode-2}
& A \, p'' + 2 \, A' \, p' + \frac{1}{s} A \, p'  +   (\Omega - \beta  \, A^2) \,  A  = 0.
\end{align} 
Here the rotation frequency $\Omega \in \mathbb{R}$ is an unknown quantity that we have to determine.

The approach developed in \cite{Gr80, Ha82, KoHo81} for solving (\ref{ode-1}--\ref{ode-2}) involves two steps.

\paragraph{Step 1: shooting arguments} Begin with the case $\beta = 0$. It follows that $\Omega =0$ and $p'(s)$ is identically zero are necessary for the existence of bounded solutions with the prescribed condition $p'(0) = 0$; see \cite{Ha82} Section 3. Hence it suffices to solve 
\begin{equation} \label{eq-a}
A'' + \frac{1}{s} A' - \frac{m^2}{s^2} A + (1- A^2) \, A = 0.
\end{equation}
The location $s = 0$ of the tip is a singularity of (\ref{eq-a}). By analyticity, any bounded nontrivial solution $A(s)$ of (\ref{eq-a}) admits the asymptotic expansion,
\begin{equation}
A(s) = a_m  s^{m} + O(s^{m+1}) \quad \mbox{as   } s \searrow 0.
\end{equation}
This solution is shown to continue globally, with $a_m \neq 0$ as the shooting parameter.

\paragraph{Step 2: perturbation arguments} The cases $0<|\beta| \ll 1$ can be treated by careful phase portrait analysis. In particular, it is shown that the shooting manifold and the center-unstable manifold of the trivial solution intersect transversely; see \cite{KoHo81} Theorem 3.1. Therefore, their intersection, which produces $m$-armed spirals, persists for small perturbations $-i \, \beta \, |\Psi|^2 \, \Psi$ on the reaction kinetic of (\ref{gle-1}).

Motivated by the existence of spiral waves on $\mathbb{R}^2$, in \cite{Paetal94} the authors considered (\ref{gle-1}) on the unit disk $\mathcal{B}^2$ equipped with Neumann boundary conditions,
\begin{equation} \label{gle3}
\partial_t \Psi  = \frac{1}{\lambda} \Delta_{\mathcal{B}^2} \Psi + ( 1 - |\Psi|^2 - i\, \beta \, |\Psi|^2 ) \, \Psi,
\end{equation}
and added a diffusion parameter $\lambda > 0$, which can also be interpreted as a scale on the radius of the disk. Based on numerical evidences, they conjectured that $m$-armed spiral waves exist for all 
\begin{equation}
\lambda > j_{0, m}^2.
\end{equation}
Here $j_{0, m}$ is the first positive zero of the derivative of the Bessel function. This conjecture was confirmed rigorously by shooting arguments in \cite{Ts10}.

There are three reasons for studying spiral waves on compact two-dimensional manifolds (e.g., the unit disk $\mathcal{B}^2$). First, in experiments and numerical simulation the domain is bounded. Second, the domain size and presence of boundary may affect the existence and the pattern of spiral waves; see \cite{Bae03, GoSt03}. Third, such a study reveals dynamical changes as the domain size tends to infinity, for instance, as $\lambda \rightarrow \infty$ in (\ref{gle3}).

\subsection{Main results} 

In this paper we consider the underlying domain $\mathcal{M}$ as a compact surface of revolution in $\mathbb{R}^3$, because we require $\mathcal{M}$ to admit polar coordinates so that the spiral Ansatz is applicable. Moreover, the unit 2-sphere differs from the unit disk topologically, and by absence of boundary. Hence we distinguish two cases, either $\partial \mathcal{M}$ is nonempty (i.e., \textit{circular geometry}) or empty (i.e., \textit{spherical geometry}), to study how topological structure affects the pattern of spiral waves.

This paper is devoted to establishing a new global bifurcation approach for pattern formation of Ginzburg-Landau spiral waves. With this approach we generalize the results of existence for rigidly-rotating spiral waves, in the following four aspects:
\begin{itemize}
\item We add complex diffusion process (see $\eta \neq 0$ below in (\ref{gle-4})) and consider more general reaction kinetics.
\item We include Dirichlet boundary conditions and Robin boundary conditions in the circular geometry.
\item We obtain two new patterns: frozen spirals in the circular and spherical geometries, and 2-tip spirals in the spherical geometry.
\item In the literature only positive amplitude $A(s) > 0$ for $s > 0$ of (\ref{eq-a}) have been discussed; see \cite{Gr80, KoHo81, Paetal94, Ts10}. This restriction is not natural. Based on our bifurcation approach we are able to obtain solutions with sign-changing amplitude; see Lemma \ref{locbif} and Fig. \ref{fig3} in this paper, and also \cite{DaPh18}.
\end{itemize}

\begin{remark}
The shooting argument used in the literature \cite{Gr80, KoHo81, Paetal94, Ts10} tracks the first critical point of the amplitude $A(s)$, and therefore, such information yields monotonic amplitude and fits well with Neumann boundary conditions, only. Hence the shooting argument is not able to treat Dirichlet and Robin boundary conditions, and moreover, it cannot find solutions with sign-changing amplitude.
\end{remark}

\begin{remark}
We not only greatly extend the existence results in \cite{Paetal94, Ts10}, but also the bifurcation structure established in this paper is useful in studying other problems related to the dynamics of spiral waves, for instance, the hyperbolicity in \cite{DaPh18} and feedback control stabilization by the method in \cite{Sch16}.
\end{remark}

More precisely, consider the general complex Ginzburg-Landau equation (abbr. GLe),
\begin{equation} \label{gle-4}
\partial_t \Psi  = \frac{1}{\lambda} (1 + i \, \eta) \Delta_{\mathcal{M}} \Psi + f(|\Psi|^2, \mathbf{b})\, \Psi.
\end{equation}
Here $\Delta_{\mathcal{M}}$ is the Laplace-Beltrami operator on a compact surface of revolution $\mathcal{M}$ to be defined shortly, $\lambda > 0$ is a bifurcation parameter, $\eta \in \mathbb{R}$ is a prescribed complex diffusion parameter, and $\mathbf{b} \in \mathbb{R}^d$ collects prescribed kinetic parameters. 

The surface of revolution is defined as
\begin{equation} \label{polar}
    \mathcal{M} = \{ \left( a(s) \cos(\varphi), a(s) \sin(\varphi), \tilde{a}(s)\right) : s \in  [0, s_*], \, \varphi \in S^1 \}.    
\end{equation}
The main examples are the unit disk when $a(s) = s$ and $\tilde{a}(s) = 0$ for $s \in [0,1]$, and the unit 2-sphere when $a(s) = \sin(s)$ and $\tilde{a}(s) = \cos(s)$ for $s \in [0,\pi]$. 

In general, we consider the smoothness class of $\mathcal{M}$ to be $C^{2, \nu}$ with a fixed exponent $\nu \in (0,1)$, and thus $a(s)$ and $\tilde{a}(s)$ are $C^{2,\nu}$ functions. Without loss of generality we consider $s$ as the arc length parameter and thus 
\begin{equation} \label{arc-length}
 (a'(s))^2 + (\tilde{a}'(s))^2 = 1 \quad \mbox{for all   } s \in [0,s_*]. 
\end{equation}
We assume
\begin{equation} \label{positivity}
a(0) = 0 \mbox{ and } a(s)>0 \mbox{ for all } s \in (0, s_*).
\end{equation}
The smoothness of $\mathcal{M}$ prevents formation of a cusp at $s = 0$ and thus implies $\tilde{a}'(0) = 0$, or equivalently, $a'(0) = 1$ by (\ref{arc-length}--\ref{positivity}). Moreover, the boundary $\partial \mathcal{M}$ is empty if and only if $a(s_*) = 0$, and in this case the smoothness of $\mathcal{M}$ implies $a'(s_*) = -1$. 

\begin{figure}[htbp] \label{fig1}
  \centering
  \label{fig:function-f}\includegraphics[scale = 1]{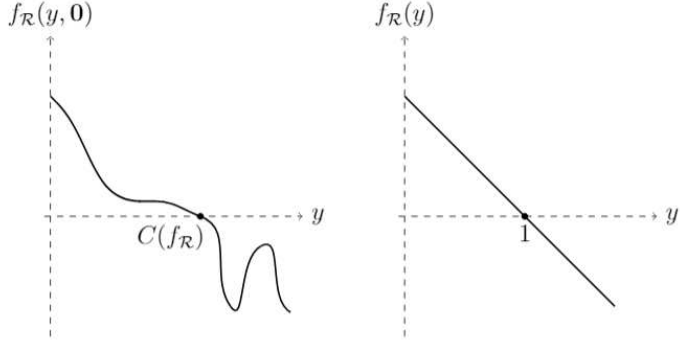}
  \caption{The graph of $f_\mathcal{R}(y,\emph{\textbf{0}})$ that satisfies (A1--A2). Left: a typical choice of $f_\mathcal{R}(y,\emph{\textbf{0}})$ for the general GLe. \textit{Right}: $f_\mathcal{R}(y) = 1-y$ for the cubic supercritical GLe.}
\end{figure}

We seek spiral wave solutions that bifurcate from the zero equilibrium $\Psi \equiv 0$. For this purpose, we consider $f \in C^3([0,\infty) \times \mathbb{R}^d, \mathbb{C})$ and impose two assumptions as we express $f = f_\mathcal{R} + i \, f_\mathcal{I}$, where $f_\mathcal{R}$ and $f_\mathcal{I}$ are real-valued functions; see Fig. \ref{fig1}.
\begin{itemize}
\item[\textbf{(A1)}] $f_{\mathcal{R}}(0, \textbf{0}) = 1$, and there exists a constant $C = C(f_\mathcal{R}) > 0$ that only depends on choices of $f_\mathcal{R}(\cdot \,, \textbf{0})$ such that 
\begin{equation}
f_\mathcal{R} (y, \textbf{0})
\left\{
\begin{array}{ll}
= 0, & y = C, \\
< 0, & y > C.
\end{array} \right.
\end{equation}
Moreover, we assume
\begin{equation} \label{fci}
f_{\mathcal{I}}(y,\textbf{0}) =0 \quad \mbox{for all   }  y \ge 0.
\end{equation}
\item[\textbf{(A2)}] $\partial_y f_\mathcal{R}(0, \textbf{0}) < 0$ and  $\partial_y f_\mathcal{R}(y, \textbf{0}) \le 0$ for all $y \in (0, C)$.
\end{itemize}

These assumptions include the prominent cubic supercritical nonlinearity,
\begin{equation}
f(y, \mathbf{b}) = 1 - y - i \, \beta \, y, \quad \mathbf{b} = \beta \in \mathbb{R},
\end{equation}
as $f_\mathcal{R}(y, \beta) = f_{\mathcal{R}}(y) = 1- y$ and $f_\mathcal{I}(y,\beta) = - \beta \, y$. Note that $C = 1$.

Indeed, the assumption (A1) yields the local instability of the zero equilibrium for the real GLe (when $\eta = 0$ and $\mathbf{b} = \mathbf{0}$), and also an a priori $C^0$-bound for solutions. We use such a $C^0$-bound to obtain global bifurcation curves of nontrivial solutions. The assumption (A2) allows us to solve the complex GLe by perturbation arguments.

We take a slightly different view on the spiral Ansatz (\ref{spians-1}),
\begin{equation} \label{spians-3}
\Psi(t,s, \varphi) = e^{-i\Omega t} \, \psi(s,\varphi), \quad \psi(s, \varphi) := u(s) \, e^{im\varphi},
\end{equation}
where the radial part $u(s)$ of $\psi$ is complex valued. We adopt the functional setting $\Delta_{\mathcal{M}}: D(\Delta_\mathcal{M})  \rightarrow L^2(\mathcal{M}, \mathbb{C})$. The domain $D(\Delta_\mathcal{M})$ is chosen to be $H^2(\mathcal{M}, \mathbb{C})$, and if $\partial \mathcal{M}$ is nonempty, also equipped with \textit{Robin boundary conditions},
\begin{equation}\label{RobinBC}
    \alpha_1  \Psi+\alpha_2  \nabla \Psi \cdot \textbf{n}= 0, 
\end{equation}
where $\alpha_1 \ge 0$ and $\alpha_2 \ge 0$ are not both zero. Here $\textbf{n}$ is the unit outer normal vector field on $\partial \mathcal{M}$. Note that $\alpha_1$ and $\alpha_2$ are required to share the same sign so that the real and imaginary parts of each solution do not grow at $\partial \mathcal{M}$. 

\begin{remark}
Robin boundary conditions for the GLe appear in the theory of superconductivity; see \cite{Duetal92, RiRu00}. Choices of boundary conditions may affect the existence of spiral waves, since bifurcating solutions are expanded by the associated eigenfunctions; see \cite{GoSt03}. For instance, with Robin boundary conditions on a disk the eigenfunctions of the GLe exhibit a spiral character, and thus in \cite{Goetal00} Section 4.2. the authors indicated that spiral waves also exist for Robin boundary conditions.
\end{remark}

The significance of considering (\ref{spians-3}) is that the $L^2$-subspace, 
\begin{equation}
L_m^2(\mathbb{C}) = \{ \psi \in L^2(\mathcal{M}, \mathbb{C}) : \psi(s, \varphi) = u(s) \, e^{im\varphi}, u(s) \in \mathbb{C} \},
\end{equation}
is invariant under the dynamics of the complex GLe, due to the global $S^1$-equivariance (\ref{glo-gau}). In particular, the following restriction is well defined:
\begin{equation} \label{proj}
\Delta_{m} := \Delta_{\mathcal{M}} \big|_{L_m^2(\mathbb{C})}: L_m^2(\mathbb{C})\rightarrow L_m^2(\mathbb{C}).
\end{equation}

Substituting the spiral Ansatz (\ref{spians-3}) into the complex GLe (\ref{gle-4}) yields the following elliptic equation for $\psi$:
\begin{equation} \label{1-fulleq}
(1 + i \, \eta) \, \Delta_{m}\psi + i \,\lambda \, \Omega \, \psi + \lambda \, f(|\psi|^2, \mathbf{b})\, \psi = 0.
\end{equation}
Here the parameters $\eta \in \mathbb{R}$, $\lambda > 0$, and $\mathbf{b} \in \mathbb{R}^d$ are given and we have to determine $\Omega \in \mathbb{R}$. We immediately see an advantage of our functional approach: Compared to the ODEs (\ref{ode-1}--\ref{ode-2}), (\ref{1-fulleq}) is an elliptic equation with bounded coefficients. 

Once a nontrivial solution pair $(\Omega, \psi)$ of (\ref{1-fulleq}) exists, we can express $u(s) = A(s) \, e^{i p(s)}$ in the polar form and thus 
\begin{equation}
\Psi(t, s,\varphi) = e^{-i \Omega t}  \left( u(s) \, e^{i m \varphi} \right) = A(s) \, e^{i (-\Omega t + p(s) + m \varphi)}.
\end{equation}
We define the \textit{pattern exhibited by $\Psi$} as the level set of zero imaginary part of its phase field, and thus given by the relation $- \Omega t + p(s) + m \varphi = 0$ (mod $\pi$). The $2\pi$-periodicity of the azimuthal angle $\varphi$ yields the following relation:
\begin{equation} \label{spiralpattern2}
\varphi =  \varphi_k(t, s) =  \frac{\Omega t - p(s)  + k \pi}{m} \, \mbox{    } \, (\mathrm{mod} \mbox{   } 2\pi), \quad \mbox{for   } k = 0,1,..., 2m-1.
\end{equation}
Then we are able to obtain the spatio-temporal pattern via polar coordinates:
\begin{equation} \label{pattern}
(t, s) \quad \mapsto \quad \bigcup_{k =0}^{2m-1} \big( a(s) \cos ( \varphi_k(t, s) ), a(s) \sin ( \varphi_k(t, s) ), \tilde{a}(s) \big).
\end{equation}
Hence we define that the solution pair \textit{exhibits a spiral} if the phase derivative $p'(s)$ is not identically zero, because the image in (\ref{pattern}) is a union of twisted curves. Otherwise it \textit{exhibits a vortex}. Clearly, a pattern is \textit{rotating} if $\Omega \neq 0$, or \textit{frozen} if $\Omega = 0$.

\begin{figure}[htbp] \label{fig-pattern}
\centering
\includegraphics[scale = 1]{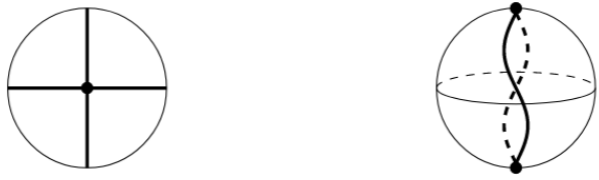}
\caption{Left: a $2$-armed vortex pattern on the disk with the origin as the vortex. Right: a $1$-armed spiral pattern on the sphere with the poles as the vortices. Both patterns may rotate rigidly.}
\end{figure}

When $\partial \mathcal{M}$ is empty, the sum of winding number associated to a tip must be zero; see \cite{Daal04}. Therefore, the topology of $\mathcal{M}$ affects the global shape of spiral patterns, and the spherical geometry supports 2-tip spirals.

\begin{theorem}[general complex GLe] \label{thm-gen}
Let $\lambda_0^{m} > 0$ be the principal eigenvalue of $-\Delta_{m}$. Then for each $\lambda > \lambda_0^{m}$, there exists an $\epsilon  > 0$ such that \emph{(\ref{1-fulleq})} possesses nontrivial solution pairs $(\Omega(\eta,\mathbf{b}), \psi(\eta,\mathbf{b}))$ parametrized by all $(\eta, \mathbf{b}) \in \mathbb{R} \times \mathbb{R}^d$ and $0 \le |\eta|, \, |\mathbf{b}| < \epsilon$. Moreover, the following statements hold:
\begin{itemize}
\item[(i)] $(\Omega(\eta,\mathbf{b}), \psi(\eta,\mathbf{b}))$ exhibits a rotating spiral if $(\eta, \mathbf{b}) \neq (0, \emph{\textbf{0}})$ lies in a small cone with $(0, \emph{\textbf{0}})$ as the vertex and $\{ (\eta, \emph{\textbf{0})}: 0 \le |\eta| < \epsilon\}$ as the axis. 
\item[(ii)] Suppose $\mathbf{b} = \beta \in \mathbb{R}$. Assume in addition
\begin{equation}
\partial_{\beta} f_\mathcal{I}(y,0) \neq 0 \quad \mbox{for all } y \in (0,C)
\end{equation}
and 
\begin{equation}
f_\mathcal{I}(0,\beta) = 0 \quad \mbox{for all } \beta \in \mathbb{R}.    
\end{equation}
Then $(\Omega(\eta,\beta), \psi(\eta,\beta))$ exhibits a rotating spiral if $(\eta, \beta) \neq (0, 0)$ lies in a small cone with $(0, 0)$ as the vertex and $\{ (0, \beta): 0 \le |\beta| < \epsilon\}$ as the axis.
\end{itemize}
\end{theorem}

We emphasize that the constant $\epsilon > 0$ depends neither on $\eta \in \mathbb{R}$ nor on $\mathbf{b} \in \mathbb{R}^{d}$. For an estimate on its best lower bound, see \cite{Ts10} for the case $\mathcal{M}$ being the disk equipped with Neumann boundary conditions.

\begin{theorem}[cubic supercritical complex GLe] \label{thm-cgle}
Suppose $\mathbf{b} = \beta \in \mathbb{R}$ and $f(y,\beta) = 1- y -i \, \beta \, y$. Then there is a smooth strictly decreasing function $\widetilde{\eta} = \widetilde{\eta}(\beta)$, $\widetilde{\eta}: (-\epsilon, \epsilon) \rightarrow \mathbb{R}$ with $\widetilde{\eta}(0) = 0$ such that $(\Omega(\eta,\beta), \psi(\eta,\beta))$ exhibits a rotating spiral if $\eta \neq \beta$ and $\eta \neq \widetilde{\eta}(\beta)$, or a frozen spiral if $\eta \neq \beta$ and $\eta = \widetilde{\eta}(\beta)$. 
\end{theorem}

It is worth noting that the existence of frozen spirals requires both $\eta \neq 0$ and $\beta \neq 0$. In the literature frozen spirals are missing because the complex diffusion parameter $\eta$ was not considered. 

\begin{figure}[htbp]
  \centering
  \label{fig:parameter}\includegraphics[scale = 1]{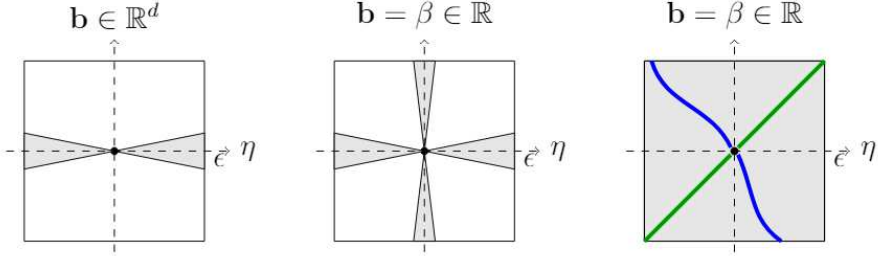}
  \caption{Types of pattern in the $(\eta,\mathbf{b})$-parameter space. The origin $(\eta, \mathbf{b}) = (0, \emph{\textbf{0}})$ supports a frozen vortex. Each parameter in the shaded region supports a rotating spiral. Left: Theorem \ref{thm-gen} (i), for which we only assume (A1--A2). Middle: Theorem \ref{thm-gen} (ii), for which we need the additional assumptions on $f_\mathcal{I}$. Right: Theorem \ref{thm-cgle}. We can completely classify the types of pattern for the cubic supercritical GLe. Each parameter on the diagonal line $\eta = \beta$ supports a rotating vortex. Each parameter on the curve $\eta = \widetilde{\eta}(\beta)$ supports a frozen spiral.}
\end{figure}

\subsection{Idea of proof}

The proof based on our functional approach
consists of three steps: global bifurcation analysis, perturbation arguments, and determination on the types of pattern.

\paragraph{Step 1: global bifurcation analysis} Begin with the unperturbed case $\eta = 0$ and $\mathbf{b} = \mathbf{0}$. We show that $\Omega = 0$ and $u(s)$ is real valued, necessarily, for every solution $\Psi$ of the spiral Ansatz (\ref{spians-3}). Hence $\Psi$ must exhibit a frozen vortex, and it suffices to solve the scalar equation,
\begin{equation} \label{1-redeq}
 \Delta_{m}\psi + \lambda \, f_\mathcal{R}(|\psi|^2, \mathbf{0}) \, \psi = 0.
\end{equation}

We next show that the spectrum of $-\Delta_{m}$ consists of simple eigenvalues, which can be listed as follows:
\begin{equation}
0 < \lambda_0^{m} < \lambda_{1}^{m} <... < \lambda_{n}^{m} < ..., \quad \lim_{n \rightarrow \infty} \lambda_{n}^{m} = \infty,
\end{equation}
and each associated eigenfunction $\mathbf{e}_n^{m}(s,\varphi) = u_n^{m}(s) \, e^{im\varphi}$ satisfies a nodal property: Its radial part $u_n^{m}(s)$ possesses exactly $n$ simple zeros in $(0,s_*)$. Therefore, nontrivial solutions of (\ref{1-redeq}) near the bifurcation point $\lambda = \lambda_n^{m} > 0$ form a unique local bifurcation curve $\mathcal{C}_n^{m}$, for each $n \in \mathbb{N}_0$; see \cite{CrRa71}. Moreover, the radial part of any nontrivial solution on $\mathcal{C}_n^{m}$ possesses exactly $n$ simple zeros. In particular, on the \textit{principal bifurcation curve} $\mathcal{C}_0^m$ the radial part does not change its sign, which allows us to extend $\mathcal{C}_0^m$ globally, in the sense that it exists for all diffusion parameter $\lambda > \lambda_0^m$, by using comparison of principal eigenvalues.

\begin{remark}
Based on the setting and approach established in \cite{Da17}, we successfully extend all other bifurcation curves globally by careful analysis on shooting manifolds; see \cite{DaPh18}. 
\end{remark}

\paragraph{Step 2: perturbation arguments} We prove that $\mathcal{C}_0^{m}$ persists for all small parameters $0 < |\eta|,\,|\mathbf{b}| \ll 1$ by applying the equivariant implicit function theorem in \cite{RePe98}.

\paragraph{Step 3: determination on the types of pattern}

We classify the types of pattern for solutions obtained in Step 2 by the \textit{frequency-parameter relation},
\begin{equation} \label{fre-par-rel}
\int_\mathcal{M} \left( \Omega - \eta \, f_\mathcal{R}(|\psi|^2, \mathbf{b}) +  f_\mathcal{I}(|\psi|^2,  \mathbf{b}) \right) |\psi|^2 =0.
\end{equation}
To derive this relation (\ref{fre-par-rel}), we multiply (\ref{1-fulleq}) by the complex conjugate $\overline{\psi}$, integrate over $\mathcal{M}$, note that $\nabla \psi \cdot \textbf{n} = u'(s) \, e^{im \varphi}$ and thus Robin boundary conditions (\ref{RobinBC}) read $\alpha_1 \, u(s) + \alpha_2 \, u'(s) = 0$, which implies that $(\nabla \psi \cdot \textbf{n} ) \, \overline{\psi} = u'(s) \, \overline{u}(s)$ is real valued, and finally sort out the real and imaginary parts.

The remaining part of this paper is organized as follows. In Section \ref{sec;2}, we establish our functional setting and derive several a priori properties of nontrivial solutions. In Section \ref{sec;3}, for the unperturbed case $\eta = 0$ and $\mathbf{b} = \mathbf{0}$, we obtain a global bifurcation curve of frozen vortices. In Section \ref{sec;4}, we show that the bifurcation curve persists for the perturbed cases $0 < |\eta|,\,|\mathbf{b}| \ll 1$. Then in Section \ref{sec;5}, we classify the parameter regime that supports spirals. Lastly, we discuss the features of spirals and indicate two problems in Section \ref{sec;6}.

\section{Preparation} \label{sec;2}

For each fixed $m \in \mathbb{N}$, we rewrite the elliptic equation (\ref{1-fulleq}) as the functional equation,
\begin{equation} \label{fulleq}
\mathcal{F}(\lambda | \, \Omega, \psi ; \, \eta, \mathbf{b}) := (1 + i \, \eta) \, \Delta_{m}\psi + i \,\lambda \, \Omega \, \psi + \lambda \, f(|\psi|^2,  \mathbf{b}) \, \psi = 0.
\end{equation}
Notice that by the global $S^1$-equivariance (\ref{glo-gau}) and the invariance (\ref{proj}), the domain and range of $\mathcal{F}$ are given by
\begin{equation}
\mathcal{F} : (0,\infty) \times \mathbb{R} \times H^2_{m}(\mathbb{C}) \times \mathbb{R} \times \mathbb{R}^d \rightarrow L_m^2(\mathbb{C}),
\end{equation}
where $H_m^2(\mathbb{C}) = L_m^2(\mathbb{C}) \cap H^2(\mathcal{M}, \mathbb{C})$ is an $H^2$-subspace.

Here we use the notation $(\lambda|\, \Omega, \psi; \, \eta, \mathbf{b})$ to distinguish $\lambda > 0$ as the bifurcation parameter; $\Omega \in \mathbb{R}$ and $\psi \in H_{m}^2(\mathbb{C})$ as the unknowns of (\ref{fulleq}); $\eta \in \mathbb{R}$ and $\mathbf{b} \in \mathbb{R}^d$ as the parameters of the complex GLe.  

Recall that $\mathcal{M}$ is a manifold of the smoothness class $C^{2, \nu}$ with a fixed $\nu \in (0,1)$.

\begin{lemma}[a priori regularity] \label{lem-apr}
Suppose that $\psi \in H_{m}^2(\mathbb{C})$ solves \emph{(\ref{fulleq})}. Then $\psi \in C^{2,\nu}(\mathcal{M},\mathbb{C})$. 
\end{lemma}
\begin{proof}
Since $H^2(\mathcal{M},\mathbb{C})$ is continuously embedded into $C^{0, \nu}(\mathcal{M}, \mathbb{C})$, from (\ref{fulleq}) we see $\psi \in C^{0, \nu}(\mathcal{M}, \mathbb{C})$. Hence $\psi \in C^{2, \nu}(\mathcal{M}, \mathbb{C})$ by using a partition of unity and the elliptic regularity theory.
\end{proof}

Zeros of nontrivial solutions are crucial in our bifurcation analysis. In polar coordinates (\ref{polar}) by the chain rule we obtain
\begin{equation}
\Delta_{\mathcal{M}} = \partial_{ss} + \frac{a'(s)}{a(s)} \, \partial_s + \frac{1}{a^2(s)} \, \partial_{\varphi \varphi},
\end{equation}
and thus $\Delta_m = \Delta_{\mathcal{M}}|_{L_m^2(\mathbb{C})}$ acts as follows:
\begin{equation} \label{proj-operator}
\Delta_m \left(u(s) \, e^{im \varphi}\right) = \left( u''(s) + \frac{a'(s)}{a(s)} u'(s) - \frac{m^2}{a^2(s)} u(s) \right) e^{i m \varphi},
\end{equation}
and thus possesses unbounded coefficients at $s = 0$, and also at $s = s_*$ if $\partial \mathcal{M}$ is empty. Therefore, we distinguish three kinds of zeros:
\begin{itemize}
\item zeros at tips $s = 0$, and also $s = s_*$ if $\partial \mathcal{M}$ is empty;
\item zeros in the interior $(0, s_*)$;
\item zeros at the boundary $s = s_*$, when $\partial \mathcal{M}$ is nonempty.
\end{itemize}

We show that zeros at tips must be isolated, and other kinds of zeros, if they exist, must be simple. The proof is based on phase portrait analysis near tips and the uniqueness of ODE initial value problems.

\begin{lemma}[characterization of zeros] \label{asym}
Suppose that $\psi \in C^{2, \nu}(\mathcal{M}, \mathbb{C})$, $\psi(s,\varphi) = u(s) \, e^{im\varphi}$, is a nontrivial solution of \emph{(\ref{fulleq})}. Then 
\begin{itemize}
\item[(i)] $s = 0$ is an isolated zero of the real and imaginary parts of $u(s)$. If in addition $\partial \mathcal{M}$ is empty, then $s = s_*$ is an isolated zero of the real and imaginary parts of $u(s)$.
\item[(ii)] $s \in (0, s_*)$ cannot be a multiple zero of $u(s)$. 
\item[(iii)] Let $\partial \mathcal{M}$ be nonempty. Then $s = s_*$ cannot be a multiple zero of $u(s)$.
\end{itemize}
\end{lemma}
\begin{proof}
Let $u = u_{\mathcal{R}} + i \, u_{\mathcal{I}}$ where $u_\mathcal{R}$ and $u_{\mathcal{I}}$ are real-valued functions. Moreover, let
$\psi= \psi_\mathcal{R} + i \, \psi_\mathcal{I}$ where $\psi_\mathcal{R}(s,\varphi) = u_\mathcal{R}(s) \, e^{im\varphi}$ and $\psi_\mathcal{I}(s,\varphi) = u_\mathcal{I}(s) \, e^{im\varphi}$. We rewrite (\ref{fulleq}) as follows:
\begin{align} \label{fulleq-2}
\begin{split}
\Delta_{m} \psi_\mathcal{R} - \eta \, \Delta_{m}\psi_\mathcal{I} + \lambda \left( -\Omega \, \psi_\mathcal{I} + f_\mathcal{R}(|\psi|^2,\mathbf{b}) \, \psi_\mathcal{R} - f_\mathcal{I}(|\psi|^2,\mathbf{b}) \, \psi_\mathcal{I} \right) &= 0,
\\
\eta \, \Delta_{m} \psi_\mathcal{R} + \Delta_{m}\psi_\mathcal{I} + \lambda \left( \Omega \, \psi_\mathcal{R} + f_\mathcal{I}(|\psi|^2,\mathbf{b}) \, \psi_\mathcal{R} + f_\mathcal{R}(|\psi|^2,\mathbf{b}) \, \psi_\mathcal{I} \right) &= 0.
\end{split}
\end{align}
Here $|\psi|^2 = u_{\mathcal{R}}^2 + u_{\mathcal{I}}^2$. The Euler multiplier 
\begin{equation} \label{euler-mul}
\left(\frac{\mathrm{d} s}{\mathrm{d} \tau} = \right) \, \, \dot{s} = a(s),
\end{equation}
transforms (\ref{fulleq-2}) to 
\begin{align}
\begin{split}
\ddot{u_{\mathcal{R}}} - m^2 \, u_\mathcal{R} - \eta \, (\ddot{u_\mathcal{I}} - m^2 \, u_\mathcal{I}) + \lambda \, a^2 \left( -\Omega \, u_\mathcal{I} + f_\mathcal{R}(|\psi|^2,\mathbf{b}) \, u_\mathcal{R} - f_\mathcal{I}(|\psi|^2,\mathbf{b}) \, u_\mathcal{I} \right) &= 0, 
\\
\eta\, (\ddot{u_{\mathcal{R}}} - m^2 \, u_\mathcal{R}) + \ddot{u_\mathcal{I}} - m^2 \, u_\mathcal{I} + \lambda \, a^2 \left( \Omega \, u_\mathcal{R} + f_\mathcal{I}(|\psi|^2,\mathbf{b}) \, u_\mathcal{R} + f_\mathcal{R}(|\psi|^2,\mathbf{b}) \, u_\mathcal{I} \right) &= 0, 
\end{split}
\end{align}

and thus the following ODE system for $\tau \in \mathbb{R}$:
\begin{align} \label{fulleq-100}
\begin{split}
\dot{u_\mathcal{R}} &= v_\mathcal{R},
\\
 \dot{v_\mathcal{R}} & =  m^2 \, u_\mathcal{R} - \lambda \, a^2(s) ( G_1 + \eta \, G_2),
\\
\dot{u_\mathcal{I}} &= v_\mathcal{I},
\\
 \dot{v_\mathcal{I}} & =  m^2 \, u_\mathcal{I} - \lambda \, a^2(s) ( - \eta \, G_1 + G_2),
\\
\dot{s} &= a(s),
\end{split}
\end{align}
where $G_1$ and $G_2$ are given by
\begin{align}
\begin{split}
G_1 = G_1(u_\mathcal{R}, u_\mathcal{I}, \Omega, \eta,\mathbf{b}) &= \frac{1}{1+\eta^2} \left( -\Omega \, u_\mathcal{I} + f_\mathcal{R}(|\psi|^2, \mathbf{b}) \, u_\mathcal{R} - f_\mathcal{I}(|\psi|^2, \mathbf{b}) \, u_\mathcal{I}\right),
\\
G_2 = G_2(u_\mathcal{R}, u_\mathcal{I}, \Omega, \eta,\mathbf{b}) & = \frac{1}{1+\eta^2} \left(\Omega \, u_\mathcal{R} + f_\mathcal{I}(|\psi|^2, \mathbf{b}) \, u_\mathcal{R} + f_\mathcal{R}(|\psi|^2, \mathbf{b}) \, u_\mathcal{I} \right).
\end{split}
\end{align}

For (i), note that $\tau = - \infty$ corresponds to $s =0$ by the Euler multiplier (\ref{euler-mul}). Since $\psi_{\mathcal{R}}(\tau, \varphi) = u_{\mathcal{R}}(\tau) \, e^{im\varphi}$ and $\psi_{\mathcal{I}}(\tau, \varphi) = u_{\mathcal{I}}(\tau) \, e^{im\varphi}$ are continuous at $\tau = - \infty$, we obtain
\begin{equation} \label{cont-at-0}
\lim_{\tau \rightarrow -\infty} u_\mathcal{R}(\tau) = \lim_{\tau \rightarrow -\infty}u_\mathcal{I}(\tau) = 0.
\end{equation}

We next claim 
\begin{equation} \label{deri-at-0}
\lim_{\tau \rightarrow -\infty} \dot{u_\mathcal{R}}(\tau) = \lim_{\tau \rightarrow -\infty} \dot{u_\mathcal{I}}(\tau) = 0.
\end{equation}
It suffices to prove the case for $u_{\mathcal{R}}$, since the proof for the other case $u_{\mathcal{I}}$ is similar. Since $\psi \in C^{2,\nu}(\mathcal{M},\mathbb{C})$, and thus in particular,
\begin{equation} \label{first}
|\nabla\psi|_{C^0} = \sup_{s \in [0,s_*]} \left( |u'(s)|^2 + \frac{m^2}{a^2(s)} \, |u(s)|^2 \right) < \infty,  
\end{equation}
we see that $|u'(0)|^2 = (u_\mathcal{R}'(0))^2 + (u_\mathcal{I}'(0))^2 < \infty$. Since $\dot{u_\mathcal{R}}(\tau) = u_\mathcal{R}'(s) \, a(s)$ by the chain rule, the conditions $a(0) = 0$ and $|u_{\mathcal{R}}'(0)| < \infty$ imply
$\lim_{\tau \rightarrow -\infty}\dot{u_\mathcal{R}}(\tau) = \lim_{s \searrow 0} u_\mathcal{R}'(s) \, a(s) = 0$.

Observe that $G_j(0,0, \Omega, \eta, \mathbf{b}) =0$ holds for all $j = 1,2$, $\Omega \in \mathbb{R}$, $\eta \in \mathbb{R}$, and $\mathbf{b} \in \mathbb{R}^d$. Thus (\ref{fulleq-100}), (\ref{cont-at-0}), and (\ref{deri-at-0}) imply that as $\tau \rightarrow -\infty$ the phase portrait of (\ref{fulleq-100}) converges to the equilibrium $ (0,0,0,0,0)$. It is straightforward to show that the equilibrium is hyperbolic and the expanding direction at $\tau = -\infty$ is given by $(1, m, 1, m, 1)$. So $\tau = -\infty$ is an isolated zero of both $u_\mathcal{R}(\tau)$ and $u_\mathcal{I}(\tau)$, or equivalently, $s = 0$ is an isolated zero of both $u_\mathcal{R}(s)$ and $u_\mathcal{I}(s)$.

If $\partial \mathcal{M}$ is empty, then similar arguments show that the phase portrait of (\ref{fulleq-100}) converges to the hyperbolic equilibrium $(0,0,0,0,s_*)$ as $\tau \rightarrow \infty$, and follows the contracting direction $(1,-m, 1, -m, -1)$. This proves (i).
 
For (ii), suppose the contrary that $\tilde{s} \in (0, s_*)$ were a multiple zero of $u(s)$, that is, $u(\tilde{s}) = u'(\tilde{s}) = 0$. Let $\tilde{\tau} = \tau(\tilde{s}) \in \mathbb{R}$ be solved by the Euler multiplier. Then $u(\tilde{\tau}) = 0$ and  $\dot{u}(\tilde{\tau}) = u'(\tilde{s}) \, a(\tilde{s}) = 0$. Then the uniqueness of ODE initial value problems implies that both $u_\mathcal{R}(\tau)$ and $u_\mathcal{I}(\tau)$ are identically zero, which is a contradiction.

For (iii), suppose the contrary that $u(s_*) = u'(s_*) = 0$. Since $a(s_*) > 0$, $s = s_*$ does not yield unbounded coefficients in (\ref{fulleq-100}). Hence the uniqueness of ODE initial value problems applies, and so both $u_\mathcal{R}(\tau)$ and $u_\mathcal{I}(\tau)$ are identically zero, which is a contradiction. 
\end{proof}

Consider a nontrivial solution $\psi \in C^{2,\nu}(\mathcal{M},\mathbb{C})$, $\psi(s,\varphi) = u(s) \, e^{im\varphi}$, of (\ref{fulleq}). Let $s \in [0, s_*]$ be a zero of $u(s)$. By Lemma \ref{asym} such a zero is isolated, and thus yields a constant $\delta > 0$ such that $\psi$ admits
the polar form,
\begin{equation} \label{polar-form-1}
\psi(s,\varphi) = \left( A(s) \, e^{ip(s)} \right) e^{im\varphi}, \quad \mbox{for all   } s \in (\tilde{s}, \tilde{s} +\delta)
\end{equation}
(or $(\tilde{s}-\delta, \tilde{s})$ if $\tilde{s} = s_*$). By substituting the polar form (\ref{polar-form-1}) into (\ref{fulleq}), using the formula (\ref{proj-operator}), and sorting out the resulting real and imaginary parts, we obtain
\begin{align} \label{full-sys-1}
& A'' + \frac{a'}{a}  A' - \frac{m^2}{a^2} A -A (p')^2 - \eta  \left( A  p'' + 2  A'  p' + \frac{a'}{a} A  p' \right) + \lambda  f_{\mathcal{R}} (|A|^2,\mathbf{b}) A  = 0,
\\ \label{full-sys-2}
&  A p'' + 2 A'  p' + \frac{a'}{a} A p' + \eta  \left(A'' + \frac{a'}{a}  A' - \frac{m^2}{a^2} A -A (p')^2 \right) + \lambda  (\Omega + f_\mathcal{I}(|A|^2, \mathbf{b}) )  A  = 0.
\end{align} 
Combining (\ref{full-sys-1}--\ref{full-sys-2}) yields
\begin{equation} \label{full-sys-3}
(1 + \eta^2) \left( A \, p'' + 2 \, A' \, p' + \frac{a'}{a} A \, p' \right) = - \lambda  \left( \Omega - \eta \, f_\mathcal{R}(|A|^2, \mathbf{b}) + f_\mathcal{I}(|A|^2, \mathbf{b})\right) A.
\end{equation}

Multiplying (\ref{full-sys-3}) by $a(s)A(s)$ and using the identity,
\begin{equation}
(a \, A^2 \, p')' = a\, A \left( A \, p'' + 2 \, A' \, p' + \frac{a'}{a} A \, p' \right),
\end{equation}
we obtain
\begin{equation} \label{trick}
(a \, A^2 \, p')' = \frac{-\lambda}{1+\eta^2} \, a\, A^2 \left( \Omega - \eta \, f_\mathcal{R}(|A|^2, \mathbf{b}) + f_\mathcal{I}(|A|^2, \mathbf{b})\right).
\end{equation}

To integrate (\ref{trick}) over $(\tilde{s},s)$ with any fixed $s \in (\tilde{s},\tilde{s}+\delta)$, we have to show that the one-sided limit $\lim_{s \searrow \tilde{s}} a(s)A^2(s) \,p'(s)$ exists. To see it, noting 
\begin{equation} \label{u'}
|u'(s)|^2 = (A'(s))^2 + ( A(s) \, p'(s) )^2, \quad \mbox{for all   } s \in (\tilde{s}, \tilde{s}+\delta),
\end{equation}
and $|u'(s)| < \infty$ for all $s \in [0,s_*]$ due to (\ref{first}), we see that $\lim_{s \searrow \tilde{s}} A(s) \, p'(s)$ exists. Hence $\lim_{s \searrow \tilde{s}}A(s) = 0$ implies
\begin{equation} \label{trick2}
\lim_{s \searrow \tilde{s}} a(s) \, A^2(s)\,p'(s) = 0.
\end{equation}
Therefore, integrating (\ref{trick}) over $(\tilde{s}, s)$ and using (\ref{trick2}) yield the following expression for the phase derivative:
\begin{equation} \label{int-rel}
p'(s) =  \frac{- \lambda}{(1+\eta^2) \, a(s) A^2(s)}
\int_{\tilde{s}}^s a(\sigma) A^2(\sigma) \left( \Omega - \eta \, f_\mathcal{R}(|A(\sigma)|^2, \mathbf{b}) + f_\mathcal{I}(|A(\sigma)|^2, \mathbf{b}) \right) \mathrm{d} \sigma.
\end{equation}

Indeed, we can extend $p'(s)$ as a continuous function on the whole interval $[0, s_*]$.

\begin{lemma} \label{lem-con}
Suppose that $\psi \in C^{2, \nu}(\mathcal{M}, \mathbb{C})$, $\psi(s,\varphi) = u(s) \, e^{im\varphi}$, is a nontrivial solution of \emph{(\ref{fulleq})}. Let $\tilde{s} \in [0,s_*]$ be a zero of $u(s)$. Then in the polar form of $\psi$ the phase derivative $p'(s)$ is continuously extendable to $s= \tilde{s}$ by the limit $
\lim_{s \rightarrow \tilde{s}}p'(s) = 0$.
\end{lemma}
\begin{proof}
Applying L'H\^{o}pital's rule on (\ref{int-rel}), we see $\lim_{s \rightarrow \tilde{s}} p'(s) =0 $ if 
\begin{equation} \label{oneside}
\lim_{s \rightarrow \tilde{s}} \frac{a(s)}{a'(s) + 2 \, a(s) \frac{A'(s)}{A(s)}} 
\end{equation}
exists and is equal to zero. 

Consider the case $\tilde{s} = 0$. Since $a'(0) = 1$ and $\lim_{s \searrow 0}A'(s) \ge 0$, the denominator in (\ref{oneside}) is positive as $s \searrow 0$. Since $a(0) = 0$, the limit (\ref{oneside}) exists and is equal to zero. The proof for the case $\tilde{s} = s_*$ is similar as we note $a'(s_*) = -1$ and $\lim_{s \nearrow s_*} A'(s) \le 0$.

Consider the case $\tilde{s} \in (0, s_*)$. Then $u(\tilde{s}) =0$ and $u'(\tilde{s}) \neq 0$ by Lemma \ref{asym} (ii). Hence $\lim_{s \searrow \tilde{s}} A(s) = 0$, and $\lim_{s \searrow \tilde{s}} A'(s) \neq 0$ by (\ref{u'}). Indeed, $\lim_{s \searrow \tilde{s}} A'(s) > 0$. Since $a(\tilde{s}) > 0$, the denominator in (\ref{oneside}) tends to the positive infinity as $s \searrow \tilde{s}$, and thus the limit (\ref{oneside}) exists and is equal to zero.
\end{proof}

We close this section by showing that the real GLe (when $\eta = 0$ and $\mathbf{b} = \mathbf{0}$) yields frozen vortices, only. This result is called the \textit{decoupling effect} because the functional equation (\ref{fulleq}) becomes essentially a scalar equation.
\begin{lemma}[decoupling effect] \label{lem-dec-eff}
Consider \emph{(\ref{fulleq})} with $\eta = 0$ and $\mathbf{b} = \emph{\textbf{0}}$. Suppose that $\psi \in C^{2, \nu}(\mathcal{M}, \mathbb{C})$, $\psi(s,\varphi) = u(s) \, e^{im\varphi}$, is a nontrivial solution. Then $\Omega = 0$, and in the polar form of $\psi$ the phase derivative $p'(s)$ is identically zero. Consequently, it suffices to seek $u(s)$ as a real-valued function, due to the global $S^1$-equivariance \emph{(\ref{glo-gau})}.
\end{lemma}
\begin{proof} 
Since $\eta = 0$ and $\mathbf{b} = \textbf{0}$, $f_\mathcal{I}$ is identically zero by the assumption (\ref{fci}), and thus $\Omega = 0$ by the frequency-parameter relation (\ref{fre-par-rel}). Then $p'(s)$ is identically zero by (\ref{int-rel}) and Lemma \ref{lem-con}.
\end{proof}

\section{Global bifurcation analysis} \label{sec;3}

In this section we consider the real GLe (when $\eta = 0$ and $\mathbf{b} = \mathbf{0}$). Due to the decoupling effect (see Lemma \ref{lem-dec-eff}), we introduce \begin{equation}
L_m^2(\mathbb{R}) := \{ \psi \in L^2( \mathcal{M}, \mathbb{C}) : \psi(s, \varphi) = u(s) \, e^{im \varphi}, u(s) \in \mathbb{R}\},
\end{equation}
and the functional equation (\ref{fulleq}) becomes
\begin{equation} \label{fulleq-3}
\mathcal{F}(\lambda | \, 0, \psi ; \, 0, \mathbf{0}) = \Delta_{m}\psi + \lambda \, f_\mathcal{R}(|\psi|^2,  \mathbf{0}) \, \psi = 0,
\end{equation}
where $\psi \in H_m^2(\mathbb{R}) := L_m^2(\mathbb{R}) \cap H^2(\mathcal{M}, \mathbb{C})$.

\begin{remark}
We identify $\mathbb{C}$ as a real vector space; the inner product of $L_m^2(\mathbb{C})$, now as a real vector space, is defined by
\begin{equation}
\langle \psi_1, \psi_2 \rangle_{L_{m}^2(\mathbb{C})} := \mathrm{Re} \, \left( \int_0^{s_*} u_1(s) \, \overline{u_2(s)} \,  a(s) \, \mathrm{d}s \right),
\end{equation}
where $\psi_j(s,\varphi) = u_j(s) \, e^{im\varphi}$ for $j = 1,2$. Hence a Fr\'{e}chet derivative is a bounded linear operator over $\mathbb{R}$. We do not identify $\mathbb{C}$ as $\mathbb{R}^2$ because multiplication by complex numbers yields simplicity of notations. 
\end{remark}

To seek nontrivial solutions that bifurcate from the trivial solution $\psi \equiv 0$, we study the following partial Fr\'{e}chet derivative of $\mathcal{F}$:
\begin{equation}
 \mathcal{L}_{(\lambda|\, 0)} := D_\psi \mathcal{F}(\lambda|\, 0, 0;\, 0,\textbf{0}): H_{m}^2(\mathbb{R}) \subset L_m^2(\mathbb{R}) \rightarrow L_m^2(\mathbb{R}),  
\end{equation}
which due to $f_\mathcal{R}(0, \mathbf{0}) = 1$ is given by
\begin{equation}
\mathcal{L}_{(\lambda|\, 0)}[U] = \Delta_m U + \lambda \, U, \quad U \in H_m^2(\mathbb{R}).
\end{equation}

We observe from (\ref{proj-operator}) that $-\Delta_m$ is a singular Sturm-Liouville operator because $a(0) = 0$ (and also $a(s_*) = 0$ if $\partial \mathcal{M}$ is empty). However, it is singular merely due to polar coordinates. We show that indeed $-\Delta_m$ shares its spectral property and nodal property of eigenfunctions with regular Sturm-Liouville operators.

\begin{lemma}[spectral property] \label{lem-spec} The spectrum of $-\Delta_{m}$ consists of simple eigenvalues, which can be listed as
\begin{equation} \label{eig-list}
0 < \lambda_0^{m} < \lambda_1^{m} < \lambda_2^{m} < ... < \lambda_n^{m} < ..., \quad \lim_{n \rightarrow \infty} \lambda^{m}_{n} =  \infty.
\end{equation}
In particular, $\Delta_{m}: H^2_{m}(\mathbb{R}) \rightarrow L_m^2(\mathbb{R})$ is a linear homeomorphism.
\end{lemma}
\begin{proof}
Since $-\Delta_{m}$ is self-adjoint, nonnegative, and has compact resolvent, the spectrum of $-\Delta_{m}$ consists of nonnegative real eigenvalues, which counting the multiplicity can be listed as $0 \le \lambda_0^m \le \lambda_2^m \le ... \le \lambda_n^m \le ...$ such that $\lim_{n \rightarrow \infty} \lambda_n^m = \infty$.

We consider the eigenvalue problem for $U \in H_{m}^2(\mathbb{R})$,
\begin{equation} \label{eig-prob}
-\Delta_{m}U = \lambda \, U, \quad U (s, \varphi) = v(s) \, e^{i m \varphi}.
\end{equation}

To show that zero is not an eigenvalue of $-\Delta_{m}$, we suppose the contrary that $U$ were nontrivial and would satisfy $-\Delta_{m}U = 0$. Multiplying $-\Delta_m U = 0$ by the complex conjugate $\overline{U}$ and integrating over $\mathcal{M}$ yield $\int_{\mathcal{M}} |\nabla U|^2 = 0$. Hence $U$ is a constant function on $\mathcal{M}$. Since $U(s,\varphi) = v(s) \, e^{im\varphi}$ and $m \neq 0$, we see that $U(s,\varphi)$ is identically zero as it is a constant function, which is a contradiction.

We next prove that all eigenvalues are simple. Since $\Delta_{m}$ is self-adjoint, it suffices 
to show that its geometric multiplicity is one. 

The eigenvalue problem (\ref{eig-prob}) is equivalent to the second-order ODE
\begin{equation} \label{eig-prob-2}
v'' + \frac{a'}{a} v' - \frac{m^2}{a^2} v = -\lambda \, v.
\end{equation}

Applying the Euler multiplier $\mathrm{d}s/\mathrm{d}\tau = \dot{s} = a(s)$
on (\ref{eig-prob-2}) yields the ODE system 
\begin{align} \label{extended}
\begin{split}
\dot{v}(\tau) &= w(\tau),
\\
\dot{w}(\tau) &= m^2 \, v(\tau) - \lambda \, a^2(s(\tau)) \, v(\tau), 
\\
\dot{s}(\tau) &= a (s(\tau)),
\end{split}
\end{align}
for all $\tau \in (-\infty, \tau_*)$. Here $\tau_* = \tau(s_*)$ is solved by the Euler multiplier. Note that $\tau_* \in \mathbb{R}$ if $\partial \mathcal{M}$ is nonempty, and $\tau_* = \infty$ if otherwise.

The key idea is to study the behavior of an eigenfunction $U(\tau, \varphi) = v(\tau) \, e^{im\varphi}$ near $\tau  = -\infty$. With a very similar proof in Lemma \ref{asym} (i), it follows that the phase portrait of (\ref{extended}) converges to the hyperbolic equilibrium $(0,0,0)$ as $\tau \rightarrow -\infty$, and moreover, $(1, m)$ is the $(v, \dot{v})$-component of the expanding direction at $\tau = -\infty$. Hence there exists a constant $c \neq 0$ such that
\begin{equation} \label{rate}
\lim_{\tau \rightarrow -\infty} \frac{v(\tau)}{e^{m\tau}} = c.
\end{equation}
\indent 
Let $U_1 (\tau,\varphi) = v_1(\tau) \, e^{im \varphi}$ and $U_2(\tau,\varphi) = v_2(\tau) \, e^{im \varphi}$ be two eigenfunctions associated with the same eigenvalue. By (\ref{rate}) there exist a constant $\delta \neq 0$
such that 
\begin{equation}
\lim_{\tau \rightarrow -\infty} \frac{v_1(\tau)}{v_2(\tau)} = \delta.
\end{equation}
Define $z(\tau) = v_1(\tau) - \delta \, v_2(\tau)$. Then
\begin{equation} \label{rate2}
\lim_{\tau \rightarrow -\infty} \frac{z(\tau)}{e^{m\tau}} = 0.
\end{equation}
We complete the proof by showing that $z(\tau)$ is identically zero. Suppose the contrary that $z(\tau)$ were not identically zero. Then $z(\tau) \, e^{im \varphi}$ would be an eigenfunction of $-\Delta_{m}$, and thus $\lim_{\tau \rightarrow -\infty}\frac{z(\tau)}{e^{m\tau}} \neq 0$ by (\ref{rate}), which contradicts to (\ref{rate2}). 
\end{proof}

\begin{lemma}[nodal property of eigenfunctions] \label{lem-nodal} Let $\mathbf{e}_n^{m}(s, \varphi) = v^{m}_n(s) \, e^{im\varphi}$ be an eigenfunction of $-\Delta_m$ associated with the eigenvalue $\lambda_n^{m}$. Then $v_n^{m}(s)$ possesses exactly $n$ simple zeros in $(0, s_*)$. Moreover, if $\mathcal{M}$ possesses the reflection symmetry,
\begin{equation} \label{ref-symm}
a(s) = a(s_* - s) \quad \mbox{for all   } s \in [0, s_*],
\end{equation}
then 
\begin{equation}
v_n^{m}(s) = (-1)^{n}\, v_n^{m}(s_* - s) \quad \mbox{for all   } n \in \mathbb{N}_0 \mbox{   and   } s \in [0,s_*].
\end{equation}      
\end{lemma}
\begin{proof}
The uniqueness of the initial value problem (\ref{extended}) assures that all zeros of $v(\tau)$ are simple. It suffices to prove that $v_n^{m}(\tau)$ possesses exactly $n$ zeros in $(-\infty,\tau_*)$. Since there are no multiple zeros, the following Pr\"{u}fer transformation is well defined:
\begin{align} 
\begin{split}
v(\tau) &= R(\tau)  \cos( \theta(\tau) ),
\\
\dot{v}(\tau) &= R(\tau)  \sin (\theta(\tau)).
\end{split}
\end{align}
Note that zeros of $v(\tau)$ are points of intersection between the phase portrait of (\ref{extended}) and the $\dot{v}$-axis. It is straightforward to derive the following ODE system:
\begin{align} \label{prufer}
\begin{split}
\big(\dot{\theta} = \big) \,\,\, \dot{\theta}_\lambda &= - \sin^2(\theta) + \left( m^2 - \lambda \, a^2(s) \right) \cos^2(\theta),
\\
\dot{s} &= a(s).
\end{split}
\end{align}
There are two crucial properties.
\begin{itemize}
\item[(P1)] 
The angle function $\theta_\lambda(\tau)$ is strictly decreasing at points where the phase portrait intersects the $\dot{v}$-axis.
\item[(P2)] 
If $\lambda_1 < \lambda_2$ and $\lim_{\tau \rightarrow -\infty}\theta_{\lambda_1}(\tau) = \lim_{\tau \rightarrow -\infty} \theta_{\lambda_2}(\tau)$, then
$\theta_{\lambda_1}(\tau) > \theta_{\lambda_2}(\tau)$ for all $\tau \in (-\infty, \tau_*)$.
\end{itemize}

We already know that the phase portrait of (\ref{extended}) converges to the hyperbolic equilibrium $(0,0,0)$ as $\tau \rightarrow -\infty$, and $(1, m)$ is the $(v, \dot{v})$-component of the expanding direction at $\tau = \infty$. Then it holds for all $\lambda \in \mathbb{R}$ that
\begin{equation} \label{fixedangle}
\lim_{\tau \rightarrow -\infty}\tan(\theta_\lambda(\tau)) = m.
\end{equation}

Suppose that $\partial \mathcal{M}$ is nonempty. Then the radial part $v(\tau)$ of an eigenfunction is given by a point of intersection between the phase portrait of (\ref{extended}) at $\tau_* \in \mathbb{R}$ and the Robin line, 
\begin{equation}
L^{\alpha_1, \alpha_2} = \{(v, \dot{v}, s_*) : \alpha_1 \, v + \alpha_2 \, \dot{v} = 0\},
\end{equation}
due to Robin boundary conditions. Since the slope of $L^{\alpha_1, \alpha_2}$ is nonpositive, (\ref{fixedangle}) implies that as $\tau \rightarrow - \infty$ the $(v,\dot{v})$-component of the phase portrait stays in the interior of the first or the third quadrant of the phase plane. Moreover, from (\ref{prufer}) and (P1) we see that if $\lambda \le 0$, then the $(v,\dot{v})$-component of the phase portrait gets trapped in the interior of the first or the third quadrant for all $\tau \in (-\infty, \tau_*)$. By the existence of eigenfunctions, (P1), (P2), and the fact that the slope of $L^{\alpha_1, \alpha_2}$ is nonpositive, we can tune $\lambda \nearrow \infty$, and then assure that $v_n^{m}(\tau)$ possesses exactly $n$ zeros in $\tau \in (-\infty, \tau_*)$.  
 
Suppose that $\partial \mathcal{M}$ is empty. The phase portrait of (\ref{extended}) converges to the hyperbolic equilibrium $(v,\dot{v},s) = (0,0,s_*)$ as $\tau \rightarrow \infty$. At $\tau = \infty$, the $(v,\dot{v})$-component of the contracting direction is $(1, -m)$, whose slope is negative. Hence similarly we can tune $\lambda \nearrow \infty$, and then assures that $v_n^{m}(\tau)$ possesses exactly $n$ zeros in $(-\infty, \tau_*)$. 

Due to the reflection symmetry (\ref{eig-prob-2}) is unchanged as we apply the new variable $s \mapsto s_* - s$. Since all eigenvalues are simple, either $v_n^{m}(s) = v_n^{m}(s- s_*)$ or $v_n^{m}(s) = - v_n^{m}(s-s_*)$ for all $s \in [0, s_*]$. Since $v_n^{m}(s)$ possesses exactly $n$ zeros in $(0, s_*)$, we see that $n \in \mathbb{N}_0$ is even if and only if $s_*/2$ is not a zero of $v_n^{m}(s)$, and thus if and only if $v_n^{m}(s) = v_n^{m}(s_* -s)$. 
\end{proof}

\begin{remark}
Let $\mathcal{M}$ be the unit disk. The eigenvalues of $-\Delta_m$ are $\lambda_{n}^{m} = j_{n, m}^2$, where $j_{n, m}$ is the $(n+1)$-th positive zero of the equation $\alpha_1 \, J_{m}(s) + \alpha_2 \, J_{m}'(s) = 0$. Here $J_{m}$ is the Bessel function of the first kind of index $m$. The eigenfunctions are $
\mathbf{e}_n^{m}(s, \varphi) =  J_{m} (j_{n, m} \, s) \, e^{ i m \varphi}$, where $J_{m}(j_{n, m} \, s)$ is the associated Bessel function.

Let $\mathcal{M}$ be the unit 2-sphere. The eigenvalues of $-\Delta_m$ are $\lambda_n^m = (m+n)(m+n+1)$. The eigenfunctions are $\mathbf{e}_n^m(s, \varphi) = P_{m+n}^m (\cos(s)) \, e^{i m \varphi}$, where $P_{m+n}^m$ is the associated Legendre function. 
\end{remark}

\begin{lemma}[local bifurcation curves] \label{locbif}
For each $n \in \mathbb{N}_0$ nontrivial solutions of \emph{(\ref{fulleq-3})} near $(\lambda_{n}^{m}|\, 0)$ form a unique $C^2$ local bifurcation curve, 
\begin{equation}
\mathcal{C}_n^{m} := \left\{  (\lambda_n(\sigma)|\, \sigma \, \mathbf{e}_n^{m} + v_n(\sigma)): 0 \le |\sigma| \ll 1 \right\} \subset (0,\infty) \times H_{m}^2(\mathbb{R}).
\end{equation}
Here $\mathbf{e}_n^{m}$ is the $L^2$-normalized eigenfunction of $-\Delta_{m}$ associated with $\lambda_n^{m}$. Moreover, the following statements hold:
\begin{itemize}
\item[(i)] $\lambda_n(0) = \lambda_{n}^{m}$, $D_\sigma \lambda_n (0) = 0$, and $D_{\sigma}^2 \lambda_n(0) > 0$. Hence the shape of $\mathcal{C}_n^{m}$ is a supercritical pitchfork.
\item[(ii)] $v_n(0) = 0$, $D_\sigma v_n(0) = 0$, and
\begin{equation} \label{orthogonal}
\big\langle v_n(\sigma), \mathbf{e}_n^{m} \big\rangle_{L_m^2(\mathbb{R})} = 0.
\end{equation}
Consequently, $\mathcal{C}_n^{m}$ intersects $\mathbb{R} \times \{0\}$ only at the bifurcation point $(\lambda_{n}^{m}|\, 0)$.  
\item[(iii)] \emph{(}$(\mathbb{Z}_2 \times \mathbb{Z}_2)$-equivariance\emph{)} Suppose that $\partial \mathcal{M}$ is empty and possesses the reflection symmetry \emph{(\ref{ref-symm})}. Let $\psi_n(\sigma) = \sigma \, \mathbf{e}_m^n + v_n(\sigma)$. Then 
\begin{equation} \label{z_2-equiv}
\psi_n(\sigma)(s,\varphi) = (-1)^n \, \psi_n(\sigma)(s_* - s, \varphi)
\end{equation}
holds for all $0 \le |\sigma| \ll 1$, $s \in [0,s_*]$, and $\varphi \in [0,2\pi)$.
\end{itemize} 
\end{lemma}
\begin{proof}
Except the proof for $D_\sigma \lambda_n (0) = 0$ and $D^2_\sigma \lambda_n (0) > 0$, the items (i) and (ii) follow from Lemma \ref{lem-spec} and the well-known local bifurcation results from simple eigenvalues; see \cite{CrRa71} Theorem 1.7. Note that these bifurcation curves are $C^2$ because $\mathcal{F}$ is $C^3$ real Fr\'{e}chet differentiable.

Let us substitute $\psi_n(\sigma) = \sigma \, \mathbf{e}_n^{m} + v_n(\sigma)$ into (\ref{fulleq-3}), differentiate (\ref{fulleq-3}) with respect to $\sigma$, and use the self-adjointness of $-\Delta_{m}$ on $L_m^2(\mathbb{R})$. It follows that
\begin{equation}
D_\sigma \lambda_n(0) = -\frac{1}{2} \, \big\langle D_\sigma v_n(0), \mathbf{e}_n^{m} \big\rangle_{L_m^2(\mathbb{R})} = 0
\end{equation}
due to (\ref{orthogonal}). Moreover,
\begin{equation}
D_{\sigma}^2 \lambda_n(0)= -2 \, \lambda_n^{m} \, \partial_y f_\mathcal{R}(0,\textbf{0}) \int_0^{s_*} \left|\mathbf{e}_n^{m}(s) \right|^4 \, a(s) \, \mathrm{d}s > 0,
\end{equation}
since $\partial_y f_\mathcal{R}(0, \textbf{0})  < 0$ by (A2). 

To prove (iii), note that the reflection symmetry implies that (\ref{fulleq-3}) is $(\mathbb{Z}_2 \times \mathbb{Z}_2)$-equivariant under the following actions: 
\begin{equation}
(\rho_N \psi ) (s, \varphi) := \psi(s_*-s, \varphi), \quad (\rho_{D} \psi ) (s, \varphi) := -\psi(s_*-s, \varphi).
\end{equation}

We consider the action by $\Gamma_N := \{\mathrm{id}, \rho_N\} \cong \mathbb{Z}_2$ and denote by $\mathcal{F}_N(\lambda|\,\psi)$ the restriction of $\mathcal{F}(\lambda|\,0, \psi; 0, \textbf{0})$ to the $\Gamma^N$-fixed subspace; see \cite{ChLa00} Chapter 2. If $n \in \mathbb{N}_0$ is even, then $(\lambda_n^m| \, 0)$ is a bifurcation point of the equation $\mathcal{F}_N(\lambda|\,\psi) = 0$; see Lemma \ref{lem-nodal}. Then the uniqueness of local bifurcation curves implies (\ref{z_2-equiv}). The proof for the other action $\rho_D$ is analogous.
\end{proof}

We decompose each bifurcation curve as follows:
\begin{equation}
\mathcal{C}_n^{m} = \mathcal{C}_{n,+}^{m} \cup \mathcal{C}_{n,-}^{m} \cup \left\{(\lambda_n^{m}| \, 0)\right\}.
\end{equation}
Here $\mathcal{C}_{n,+}^{m}$ is the subset of $\mathcal{C}_n^{m}$ that collects all $\sigma > 0$, and $\mathcal{C}_{n,-}^{m}$ collects all $\sigma < 0$. Note that the pitchfork shape of bifurcation curves represents the obvious $\mathbb{Z}_2$-symmetry: $(\lambda|\,\psi)$ is a solution of (\ref{fulleq-3}) if and only if $(\lambda|\,$$- \psi)$ is a solution. From now on the notation $\iota$ always applies for both cases $+$ and $-$.

Since the shape of $\mathcal{C}_{n}^{m}$ is a supercritical pitchfork, taking $|\sigma|$ sufficiently small if necessary, there exists a $\delta_n > 0$ such that  $\mathcal{C}_{n,\iota}^{m}$ admits a \textit{monotone parametrization} in $\lambda \in (\lambda_{n}^{m}, \lambda_{n}^{m}+\delta_n]$, that is, there is a smooth function $\hat{\psi}_{n,\iota} : (\lambda_{n}^{m}, \lambda_{n}^{m}+\delta_n]  \rightarrow H_m^2(\mathbb{R})$ such that $(\lambda|\, \psi) \in \mathcal{C}_{n,\iota}^{m}$ 
if and only if $\lambda \in (\lambda_{n}^{m}, \lambda_{n}^{m}+\delta_n]$ and $\psi = \hat{\psi}_{n,\iota}(\lambda)$.

Our main task is to prove that $\delta_n$ extends to infinity, and thus $\mathcal{C}_{n}^{m}$ is global and undergoes no secondary bifurcations in $(0,\infty) \times H_{m}^2(\mathbb{R})$.
 
The idea of proof is based on \textit{open-closed arguments} (or anylytic induction): Show that the set of $\delta >0$ such that $\mathcal{C}_{n,\iota}^{m}$ admits a monotone parametrization in $\lambda \in (\lambda_{n}^{m}, \lambda_{n}^{m}+\delta]$ is open and closed in $(0, \infty)$. 

We need two lemmas for our open-closed arguments. First, we label each bifurcation curve by the number of simple zeros of its associated eigenfunction.
\begin{lemma}[nodal structure of bifurcation curves] \label{lem-nod-4}
Suppose that $\mathcal{C}_{n,\iota}^{m}$ admits a monotone parametrization in $\lambda \in (\lambda_{n}^{m,(\alpha)}, \lambda_{n}^{m}+\delta]$ for some $\delta > 0$. Then for each $(\lambda|\, \psi) \in \mathcal{C}_{n,\iota}^{m}$ and $\psi(s,\varphi) = u(s)\, e^{im\varphi}$, $u(s)$ possesses exactly $n$ simple zeros in $(0,s_*)$. 
\end{lemma}
\begin{proof}
Due to a priori regularity (see Lemma \ref{lem-apr}), $u(s)$ is a nontrivial solution of the following ODE:
\begin{equation} \label{2ndode}
u'' + \frac{a'}{a}  u' - \frac{m^2}{a^2} u + \lambda \, f_{\mathcal{R}} (u^2,\textbf{0}) \, u = 0.
\end{equation}

We have two observations. First, zeros of $u(s)$ are isolated by Lemma \ref{asym}, and thus there are finitely many simple zeros of $u(s)$ in $(0,s_*)$. Second, since $\mathcal{C}_{n,\iota}^{m}$ admits a monotone parametrization, the map $\lambda \mapsto u = u_\lambda(\cdot)$ is well defined for all $\lambda \in (\lambda_{n}^{m}, \lambda_{n}^{m}+\delta]$ and smooth, that is, we have the smooth dependence of solutions on the parameter $\lambda$.

These two observations imply that all nontrivial solutions $u_\lambda(s)$ of (\ref{2ndode}) with $\lambda \in (\lambda_{n}^{m}, \lambda_{n}^{m}+\delta]$ possess the same number of simple zeros in $(0,s_*)$. Recall from Lemma \ref{lem-nodal} that the radial part $u_n^{m}(s)$ of the eigenfunction $\textbf{e}_n^{m}$ possesses exactly $n$ simple zeros in $(0, s_*)$. Therefore, nontrivial solutions $u_\lambda(s)$ with $\lambda$ sufficiently near $\lambda_n^{m}$ possess exactly $n$ simple zeros in $(0,s_*)$. 
\end{proof}

We next obtain an important $C^0$-bound for solutions. 
\begin{lemma}[$C^0$-bound] \label{apriori}
Let $(\lambda| \, \psi) \in \mathcal{C}_n^{m}$ and $\psi(s,\varphi) = u(s) \, e^{im\varphi}$. Then
\begin{equation} \label{priori}
|\psi|_{C^0} = \sup_{s \in [0,s_*]} |u(s)|  \le \sqrt{C}.
\end{equation}
\end{lemma}
\begin{proof}
Suppose the contrary that $|\psi|_{C^0} > \sqrt{C}$. Then there would be an $s_M \in [0, s_*]$ such that $|u(s_M)| > \sqrt{C}$. Since $[0,s_*]$ is compact, we assume without loss of generality that $s_M$ is a global extreme point, that is, $|u(s_M)| \ge |u(s)|$ for all $s \in [0,s_*]$. We consider the case $u(s_M) > \sqrt{C}$ only, because the proof for the other case $u(s_M) < - \sqrt{C}$ is analogous.

We show $s_M \in (0, s_*)$ by dealing with the following four cases:

Case 1: $\partial \mathcal{M}$ is empty. since $u(0) = u(s_*) = 0$ by Lemma \ref{asym} (i), $s_M \in (0, s_*)$.

If $\partial \mathcal{M}$ is nonempty, then $u(0) = 0$ by Lemma \ref{asym} (i).

Case 2: Dirichlet boundary conditions. Then $u(s_*) = 0$, and so $s_M \in (0,s_*)$.  

Case 3: Robin boundary conditions with $\alpha_1 > 0$ and $\alpha_2 > 0$. Suppose the contrary that $s_M = s_*$. Robin boundary conditions imply $u'(s_*) < 0$, and thus $s = s_*$ cannot be a global extreme point of $u(s)$, which is a contradiction. 

Case 4: Neumann boundary conditions. Suppose the contrary that $s_M = s_*$. Since $u(s_*) > \sqrt{C}$ and $u'(s_*) = 0$, by continuity there would exist a constant $c> 0$ such that $u(s) > \sqrt{C}$ and
\begin{equation} \label{ineq2}
|u'(s)| < \frac{m^2}{\max_{s\in[0,s_*]} a(s)}  \sqrt{C} \quad \mbox{for all  } s \in (s_* - c, s_*).
\end{equation}
The assumption (A1) implies $f_\mathcal{R}(|u(s)|^2, \textbf{0}) < 0$ for all $s \in (s_*-c, s_*)$. Since $\lambda > 0$, from (\ref{2ndode}) the following inequalities hold for all $s \in (s_*-c, s_*)$:
\begin{align}
\begin{split}
u''(s) & \ge \frac{m^2}{a^2(s)} \, u(s) - \frac{|a'(s)|}{a(s)}\, |u'(s)|
\\& \ge \frac{m^2}{a^2(s)} \sqrt{C} - \frac{1}{a(s)}\, |u'(s)|
\\&  = \frac{1}{a^2(s)} \bigg( m^2 \sqrt{C} - a(s) \, |u'(s)| \bigg)
\\& > 0.
\end{split}
\end{align}
Note that we have used $\sup_{s \in [0,s_*]}|a'(s)| \le 1$ for the second inequality, and (\ref{ineq2}) for the last inequality. Since $u'(s_*) = 0$ and $u''(s) > 0$ for all $s \in (s_*-c, s_*)$, we see $u'(s) < 0$ for all $s \in (s_*-c, s_*)$, and thus $s = s_*$ cannot be a global extreme point of $u(s)$, which is a contradiction.

Therefore, $s_M$ lies in the interior of $[0,s_*]$, and so $u'(s_M) = 0$ and $u''(s_M) \le 0$, yielding a contradiction as we plug $s = s_M$ into (\ref{2ndode}). 
\end{proof}

We shall extend a monotone parametrization by the implicit function theorem. Thus we show that the partial Fr\'{e}chet derivative of $\mathcal{F}$ at $(\lambda|\,\psi) \in \mathcal{C}_{n,\iota}^m$ given by
\begin{equation} \label{linear-at-solution}
\mathcal{L}_{(\lambda|\,\psi)}[U]
= 
\Delta_{m}U
+ \lambda\, f_{\mathcal{R}}(|\psi|^2,\textbf{0}) \, U
+ 2\, \lambda \, \partial_y f_\mathcal{R}(|\psi|^2, \textbf{0}) \, |\psi|^2 \, U,
\end{equation}
is a linear homeomorphism. Since $\mathcal{L}_{(\lambda|\,\psi)}$ is Fredholm of index zero (see \cite{Agetal97} Theorem 2.4.1.), it suffices to show that its kernel $\mathrm{ker}\mathcal{L}_{(\lambda|\,\psi)}$ is trivial.

Our idea, which is based on comparison of principal eigenvalues, is only able to prove openness for the principal bifurcation curve $\mathcal{C}_{0}^{m}$, since it requires that the radial part of solutions does not change sign. 

\begin{lemma}[openness] \label{openness}
Suppose that $\mathcal{C}_{0,\iota}^{m}$ admits a monotone parametrization in $\lambda \in (\lambda_0^{m}, \lambda_0^{m}+ \delta]$ for some $\delta > 0$. Then there exists a $\tilde{\delta} > \delta$ such that the monotone parametrization extends to $\lambda \in (\lambda_0^{m}, \lambda_0^{m}+ \tilde{\delta})$.
\end{lemma}
\begin{proof}
Let $\lambda_0 := \lambda_0^{m}+ \delta$ and $\psi_0 := \hat{\psi}_{0,\iota}(\lambda_0)$, where $\hat{\psi}_{0,\iota}$ is a monotone parametrization of $\mathcal{C}_{0,\iota}^{m}$. By (\ref{linear-at-solution}) the linear equation $\mathcal{L}_{(\lambda_0|\,\psi_0)}[U] = 0$ for $U \in H_m^2(\mathbb{R})$ reads
\begin{equation}
\left( \Delta_{m} 
+ \lambda_0\, f_{\mathcal{R}}(|\psi_0|^2,\textbf{0})
+ 2 \, \lambda_0\, \partial_y f_\mathcal{R}(|\psi_0|^2, \textbf{0}) \, |\psi_0|^2 \right) U = 0.
\end{equation}

To show that $U(s, \varphi)$ is identically zero, we compare the principal eigenvalues of two different but related linear operators. To see it, let $\mu_1^*$ be the principal eigenvalue of the following eigenvalue problem:
\begin{equation} \label{eig-pro1}
\left(\Delta_{m} + \lambda_0 \, f_{\mathcal{R}}(|\psi_0|^2,\textbf{0}) \right) V_1 = \mu_1 \, V_1.
\end{equation}
Since $\psi_0 \in H_m^2(\mathbb{R})$ is a nontrivial solution of (\ref{eig-pro1}) for $\mu_1 = 0$, and by Lemma \ref{lem-nod-4} its radial part does not possess zeros in $(0,s_*)$, we see $\mu_1^* = 0$.

Compare $\mu_1^*$ with the principal eigenvalue $\mu_2^*$ of the following eigenvalue problem:
\begin{equation} \label{eig-pro2}
\left( \Delta_{m}  + \lambda_0 \, f_{\mathcal{R}}(|\psi_0|^2,\textbf{0}) + 2 \,\lambda_0 \, \partial_y f_\mathcal{R}(|\psi_0|^2,\textbf{0}) \, |\psi_0|^2 \right) V_2 = \mu_2 \, V_2.
\end{equation}
The $C^0$-bound (\ref{priori}) and (A2) imply $\partial_y f_\mathcal{R}(|\psi_0(s,\varphi)|^2,\textbf{0}) \le 0$ for all $s \in (0,s_*]$ and $\varphi \in [0,2\pi)$. Moreover, $s = 0$ is an isolated zero of any nontrivial solution $V_2(s, \varphi)$ of (\ref{eig-pro2}), and so $\partial_y f_\mathcal{R}(|\psi_0(0,\varphi)|^2;\textbf{0}) = \partial_y f_\mathcal{R}(0;\textbf{0}) < 0$ by (A2). Therefore,
\begin{equation}
\int_\mathcal{M} \partial_y f_\mathcal{R}(|\psi_0|^2;\textbf{0}) \, |\psi_0|^2 \, |V_2|^2 < 0.
\end{equation}
The Rayleigh quotient for principal eigenvalues implies $\mu_2^* < \mu_1^* = 0$, and so $U(s, \varphi)$ is identically zero. 
\end{proof}

\begin{remark}
Since $\mu_2^* < 0$, $\psi_0$ is a locally asymptotically stable equilibrium of the real GLe under $H_{m}^2(\mathbb{R})$-perturbations on initial conditions. It is interesting to determine the stability of $\psi_0$ under full $H^2(\mathcal{M}, \mathbb{C})$-perturbations.
\end{remark}

We prove closedness of a monotone parametrization for every bifurcation curve. As a result, a bifurcation curve is global if openness can be proved.

\begin{lemma}[closedness] \label{closed-2}
Suppose that $\mathcal{C}_{n,\iota}^{m}$ admits a monotone parametrization in $\lambda \in (\lambda_n^{m}, \lambda_n^{m}+ \tilde{\delta})$ for some $\tilde{\delta} > 0$. Then the monotone parametrization extends to $\lambda = \lambda_n^{m}+ \tilde{\delta}$. 
\end{lemma}
\begin{proof}
Let $\lambda_n := \lambda_n^{m}+ \tilde{\delta}$, and $(\lambda^j|\, \psi^j)$ be a sequence in $\mathcal{C}_{n,\iota}^{m}$ such that $\lambda^j \nearrow  \lambda_n$. We show that there exists a $\psi_n \in H_{m}^2(\mathbb{R})$ such that $\lim_{j \rightarrow \infty} \psi^j = \psi_n$ holds in $H_{m}^2(\mathbb{R})$, $(\lambda_n|\, \psi_n) \in \mathcal{C}_{n,\iota}^{m}$, and the monotone parametrization extends to $\lambda = \lambda_n$.

The $C^0$-bound (\ref{priori}) yields the uniform $L^2$-bound,
\begin{equation} \label{l2-bdd}
|\psi^j|_{L^2}^2 = \int_{\mathcal{M}}|\psi^j|^2 \le C \, \mathrm{Vol}(\mathcal{M}).
\end{equation}

Multiplying (\ref{fulleq-3}) by the complex conjugate $\overline{\psi}$, integrating over $\mathcal{M}$, applying the divergence theorem, and using Robin boundary conditions yield 
\begin{align}
\begin{split} \label{bdy-int}
\int_\mathcal{M} |\nabla \psi^j|^2 & = \frac{-\alpha_1}{\alpha_2} \int_{\partial \mathcal{M}} |\psi^j|^2 + \lambda^j \, \int_\mathcal{M} f_\mathcal{R}(|\psi^j|^2, \textbf{0}) \, |\psi^j|^2 
\\& \le \frac{\alpha_1}{\alpha_2} \, C \, \mathrm{Vol}(\partial \mathcal{M}) + \lambda^j \,C \, \mathrm{Vol}(\mathcal{M}).
\end{split}
\end{align}
In the inequality we have used $\sup_{y \in [0,C]}|f_\mathcal{R}(y,\textbf{0})| = 1$ due to (A1--A2). Note that the boundary integral in (\ref{bdy-int}) vanishes if $\alpha_2 = 0$. By (\ref{l2-bdd}--\ref{bdy-int}) we have a uniform $H^1$-bound for solutions on every compact $\lambda$-subinterval.

Since $H^1(\mathcal{M},\mathbb{C})$ is compactly embedded into $L^2(\mathcal{M},\mathbb{C})$, passing to a subsequence if necessary, there exists a $\psi_n \in L_{m}^2(\mathbb{R})$ such that $\lim_{j \rightarrow \infty}\psi^j = \psi_n$ holds in $L_{m}^2(\mathbb{R})$. Moreover, by the triangle inequality and the $C^0$-bound (\ref{priori}), it follows that $\lim_{j \rightarrow \infty} \lambda^j \, f_\mathcal{R}(|\psi^j|^2,\textbf{0}) \, \psi^j = \lambda_n\, f_\mathcal{R}(|\psi_n|^2,\textbf{0}) \, \psi_n$ holds in $L_m^2(\mathbb{R})$. Since $\Delta_{m} : H_m^2(\mathbb{R}) \rightarrow L_m^2(\mathbb{R})$ is a linear homeomorphism, we see that
\begin{equation}
\lim_{j \rightarrow \infty} \Delta_{m}^{-1} \left( \lambda^j \, f_\mathcal{R}(|\psi^j|^2,\textbf{0}) \, \psi^j \right) = \Delta_{m}^{-1} \left( \lambda_n \, f_\mathcal{R}(|\psi_n|^2,\textbf{0}) \, \psi_n \right)
\end{equation}
holds in $H_m^2(\mathbb{R})$. Hence, $\lim_{j \rightarrow \infty}\psi^j =  \psi_n$ holds in $H_m^2(\mathbb{R})$ and $(\lambda_n|\,\psi_n)$ is also a solution of (\ref{fulleq-3}). 

By Lemma \ref{locbif}, near each bifurcation point, nontrivial solutions of (\ref{fulleq-3}) form a unique local bifurcation curve whose shape is a supercritical pitchfork. Thus $\psi_n$ is not identically zero by the implicit function theorem. Hence $(\lambda_n|\, \psi_n) \in \mathcal{C}_{n,\iota}^{m}$, and the monotone parametrization extends to $\lambda = \lambda_n$ by Lemma \ref{lem-nod-4} and Lemma \ref{openness}.
\end{proof}

\begin{lemma}[globalness] \label{lem-global}
The principal bifurcation curve $\mathcal{C}_0^{m}$ is global and undergoes no secondary bifurcations in $(0,\infty) \times H_m^2(\mathbb{R})$.
\end{lemma}
\begin{proof}
The set of $\delta >0$ such that $\mathcal{C}_{0,\iota}^{m}$ admits a monotone parametrization in $\lambda \in (\lambda_{0}^{m}, \lambda_{0}^{m}+\delta]$ is nonempty by Lemma \ref{locbif}, open in $(0, \infty)$ by Lemma \ref{openness}, and closed in $(0,\infty)$ by Lemma \ref{closed-2}.
\end{proof}

\begin{figure}[htbp] \label{fig3}
  \centering
  \label{fig:bifurcatin-diagram}\includegraphics[scale = 1]{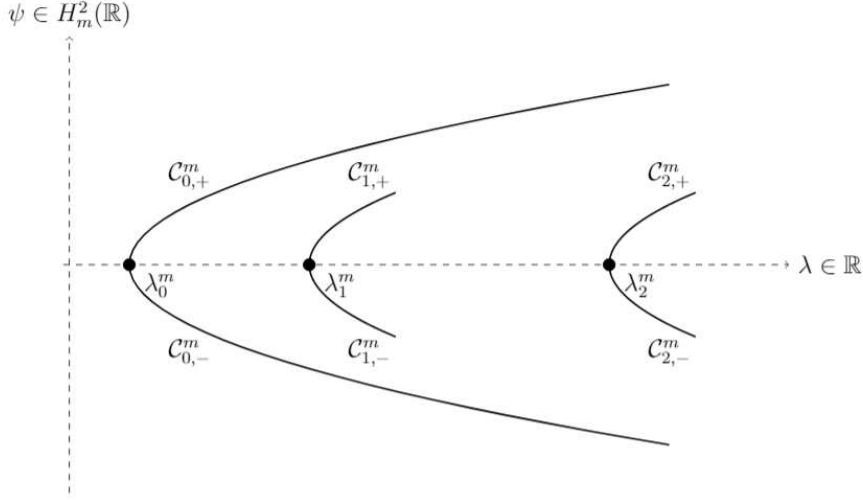}
  \caption{Bifurcation diagram $(\lambda|\,\psi)$ of the equation \emph{(\ref{fulleq-3})}. Notice that the $\psi$-component of a bifurcation curve may be unbounded, since the $H^2$-bound of solutions is not guaranteed by the $C^0$-bound \emph{(\ref{priori})}.}
\end{figure}

\section{Perturbation arguments} \label{sec;4}

For each $(\lambda_n|\, \psi_n) \in \mathcal{C}_{n, \iota}^m$ that solves 
\begin{equation}
\mathcal{F}(\lambda_n| \, 0, \psi_n; \, 0, \textbf{0}) = 0,
\end{equation}
we aim at solving the full functional equation
\begin{equation} \label{fulleq-5}
\mathcal{F}(\lambda | \, \Omega, \psi ; \, \eta, \mathbf{b}) = 0
\end{equation}
and parametrize the unknowns $\Omega \in \mathbb{R}$ and $\psi \in H_{m}^2(\mathbb{C})$ by $\lambda >0$, $\eta \in \mathbb{R}$, and $\mathbf{b} \in \mathbb{R}^d$ in the vicinity of $(\lambda|\,\eta, \mathbf{b}) =(\lambda_n|\, 0, \textbf{0})$. This task leads us to study the partial Fr\'{e}chet derivative,
\begin{equation}
\mathcal{L}_n := D_\psi \mathcal{F}(\lambda_n|\, 0, \psi_n;0, \textbf{0}) : H_{m}^2(\mathbb{C}) \subset L_m^2(\mathbb{C}) \rightarrow L_m^2(\mathbb{C}).
\end{equation}
Notice that the underlying space is $L_m^2(\mathbb{C})$, instead of $L_m^2(\mathbb{R})$, since the radial part of solutions when $(\eta, \mathbf{b}) \neq (0, \mathbf{0})$ can be complex valued. Since $\psi_n(s,\varphi) = u_n(s) \, e^{im \varphi}$ and $u_n(s)$ is real valued (see Lemma \ref{lem-dec-eff}), it is straightforward to derive
\begin{equation} \label{gat-2}
\mathcal{L}_n [U]
= 
\Delta_{m}U
+ \lambda_n\, f_{\mathcal{R}}(|\psi_n|^2,\textbf{0}) \, U
+ 2\, \lambda_n \, \partial_y f_\mathcal{R}(|\psi_n|^2, \textbf{0}) \, |\psi_n|^2 \, U_R \, e^{im\varphi}.
\end{equation}
Here $U \in H_{m}^2(\mathbb{C})$ is given by
\begin{equation}
U(s, \varphi) = (U_R(s) + i \, U_I(s)) \, e^{i m \varphi},
\end{equation}
and $U_\mathcal{R}(s)$ and $U_\mathcal{I}(s)$ are real-valued functions. It follows that  $\mathcal{L}_n$ is self-adjoint on $L_m^2(\mathbb{C})$, and moreover, Fredholm of index zero; see \cite{Agetal97} Theorem 2.4.1..

Notice that the standard implicit function theorem is not applicable for our perturbation arguments. The reason is that the global $S^1$-equivariance (\ref{glo-gau}) yields a nontrivial kernel of $\mathcal{L}_n$. Indeed, 
\begin{equation} \label{span-of-kernel}
\mathrm{span}_\mathbb{R} \langle i \, \psi_n \rangle \subset \mathrm{ker} \mathcal{L}_n. 
\end{equation}

This situation, however, is amendable to the following \textit{equivariant implicit function theorem}. The proof is essentially a restriction of the functional equation to the quotient space of a group orbit. Then the standard implicit function theorem applies.

\begin{lemma}[\cite{RePe98} Theorem 3.1] \label{eq-ift}
Let $X$, $Y$, and $M$ be Banach spaces. Suppose that $\mathcal{H} \in C^k( Y \times M, X)$ for $k \ge 2$ and $\mathcal{H}(v_0, 0) = 0$. Let $\Gamma$ be a compact Lie group and $\rho_X$ (resp. $\rho_Y$) be an action of $\Gamma$ on $X$ (resp. on $Y$). Suppose that $\mathcal{H}(\cdot, \mu)$ is $\Gamma$-equivariant, that is, $\mathcal{H}(\rho_Y(\gamma)v, \mu) = \rho_X(\gamma) \mathcal{H} (v, \mu)$ for all $\gamma \in \Gamma$, $v\in Y$, and $\mu \in M$. Assume the following:
\begin{itemize}
\item $\mathcal{L} := D_v \mathcal{H}(v_0,0)$ is Fredholm of index zero.
\item The isotropy subgroup of $v_0$ is trivial. 
\item 
Assume $M = M_1 \times M_2$ such that 
\begin{equation} \label{dim-cond-1}
\mathrm{dim} \, M_1 = \mathrm{dim} \, \mathrm{ker}  \mathcal{L} 
\end{equation}
and
\begin{equation} \label{dim-cond-2}
X = \mathcal{L} Y \oplus D_\mu \mathcal{H}(v_0,0) M_1.
\end{equation}
\end{itemize}
Then there exist a neighborhood $W \subset Y$ of the group orbit $\Gamma v_0$, neighborhoods $W_j \subset M_j$ of $\mu_j =0$ for $j = 1, 2$, and $C^k$-smooth maps $\tilde{v} : W_2 \rightarrow Y$ and $\widetilde{\mu}_1: W_2 \rightarrow M_1$ with $\tilde{v}(0) = v_0$ and $\widetilde{\mu}_1(0) = 0$ such that for each $(v, \mu_1, \mu_2) \in W \times W_1 \times  W_2$, $ \mathcal{H}(v, \mu_1, \mu_2) = 0$ if and only if $v = \rho_Y(\gamma) \tilde{v}(\mu_2)$ for some $\gamma \in \Gamma$ and $\mu_1 = \widetilde{\mu}_1(\mu_2)$.
\end{lemma}

To apply Lemma \ref{eq-ift}, we choose 
$$
\mathcal{H} := \mathcal{F}, \quad X := L_m^2(\mathbb{C}), \quad Y := H_m^2(\mathbb{C}), \quad \Omega \in M_1 := \mathbb{R}, \quad (\lambda|\, \eta, \mathbf{b}) \in M_2 := \mathbb{R}^{d+2}.
$$
Moreover, we take 
$$
v_0 := \psi_n, \quad \Gamma := S^1, \quad \mathcal{L} := \mathcal{L}_n,
$$
where $\rho_Y(\gamma) v := e^{i \gamma} v$.
It follows that the isotropy subgroup of $\psi_n$ is trivial.

We denote by $B_{\epsilon}^{d+1}(\textbf{0})$ the $(d+1)$-dimensional ball with radius $\epsilon > 0$ centered at the origin. 

\begin{lemma}[perturbation arguments] \label{lem-per-arg}
Let $(\lambda_0|\, \psi_0) \in \mathcal{C}_{0,\iota}^{m}$. Then there exists an $\epsilon_1 > 0$ and smooth functions
\begin{equation}
\widetilde{\psi} : (\lambda_0 - \epsilon_1, \lambda_0 + \epsilon_1) \times B^{d+1}_{\epsilon_1}(\emph{\textbf{0}}) \rightarrow H_{m}^2(\mathbb{C}), \quad \widetilde{\psi}(\lambda_0|\, 0, \emph{\textbf{0}}) = \psi_0,
\end{equation}
and 
\begin{equation}
\widetilde{\Omega} : (\lambda_0 - \epsilon_1, \lambda_0 + \epsilon_1) \times B^{d+1}_{\epsilon_1}(\emph{\textbf{0}}) \rightarrow \mathbb{R}, \quad  \widetilde{\Omega}(\lambda_0|\, 0, \emph{\textbf{0}}) = 0,
\end{equation}
such that for each $(\lambda|\, \eta, \mathbf{b}) \in (\lambda_0 - \epsilon_1, \lambda_0 + \epsilon_1) \times B_{\epsilon_1}^{d+1}(\emph{\textbf{0}})$, $(\widetilde{\Omega}(\lambda|\, \eta, \mathbf{b}), \widetilde{\psi}(\lambda|\, \eta, \mathbf{b}))$ is a nontrivial solution pair of \emph{(\ref{fulleq-5})}.
\end{lemma}
\begin{proof}
We apply Lemma \ref{eq-ift} and it suffices to verify (\ref{dim-cond-1}--\ref{dim-cond-2}). 

To verify (\ref{dim-cond-1}), since $M_1 = \mathbb{R}$ is real one-dimensional, and $\mathrm{span}_\mathbb{R} \langle i \, \psi_0 \rangle \subset \mathrm{ker}\mathcal{L}_0$ (see (\ref{span-of-kernel})), we need to show $\mathrm{span}_\mathbb{R} \langle i \, \psi_0 \rangle = \mathrm{ker}\mathcal{L}_0$.

By (\ref{gat-2}) the linear equation $\mathcal{L}_0 [U] = 0$ for $U \in H_{m}^2(\mathbb{C})$, $U(s, \varphi) = (U_R(s) + i\, U_I(s))\, e^{im\varphi}$, is equivalent to the following decoupled system:
\begin{align} 
& \left( \Delta_{m}
+ \lambda_0\, f_{\mathcal{R}}(|\psi_0|^2,\textbf{0})
+ 2 \,\lambda_0\,\partial_y f_\mathcal{R}(|\psi_0|^2, \textbf{0}) \, |\psi_0|^2 \right) U_R \, e^{im\varphi}  &= 0,
\\ 
&\left( \Delta_{m}
+ \lambda_0\, f_{\mathcal{R}}(|\psi_0|^2,\textbf{0}) \right) U_I \, e^{im\varphi}   & = 0.
\end{align} 
Since $(\lambda_0|\, \psi_0) \in \mathcal{C}_{0, \iota}^{m}$, the same proof in Lemma \ref{openness} shows that $U_\mathcal{R}(s) \, e^{im\varphi}$ is identically zero and $U_\mathcal{I}(s) \, e^{im\varphi} = c \, \psi_0(s,\varphi)$ for some $c \in \mathbb{R} \setminus \{0\}$. Hence $U \in \mathrm{span}_\mathbb{R} \langle i \, \psi_0 \rangle$, and so $\mathrm{ker}\mathcal{L}_0 \subset \mathrm{span}_\mathbb{R} \langle i \, \psi_0 \rangle$.

To verify (\ref{dim-cond-2}), we calculate 
\begin{equation} 
D_{(\lambda|\, \Omega; \, \eta, \mathbf{b})} \mathcal{F}(\lambda_0|\, 0, \psi_0; \, 0,\textbf{0}) (\lambda_0 |\,\Omega_*; \, 0, \textbf{0}) = i \, \lambda_0 \, \Omega_* \, \psi_0.
\end{equation}
Suppose that $V = (V_R + i\, V_I) \, e^{i m\varphi} \in H_{m}^2(\mathbb{C})$ solves $\mathcal{L}_n V = i\, \lambda_0 \, \Omega_* \, \psi_0$. Then the equation for $V_I\,e^{im\varphi}$ reads
\begin{equation}
\widetilde{L} [V_I \,e^{i m \varphi}] := \left(\Delta_{m} + \lambda_0 \, f_\mathcal{R}(|\psi_0|^2,\textbf{0}) \right) V_I \,e^{im\varphi} = \lambda_0  \, \Omega_* \, \psi_0.
\end{equation}
Since $\widetilde{L}$ is self-adjoint on $L_m^2(\mathbb{C})$ and $\widetilde{L}[\psi_0] = 0$, we have $\lambda_0 \, \Omega_* \langle \psi_0, \psi_0 \rangle_{L_m^2(\mathbb{C})} = 0$. Hence $\Omega_* = 0$, and so $(\ref{dim-cond-2})$ is verified.
\end{proof}

\section{Determination on types of pattern} \label{sec;5}

We fix $\lambda > 0$ and let 
\begin{equation}
\Omega = \widetilde{\Omega}(\eta, \mathbf{b}) \in \mathbb{R}, \quad \psi = \widetilde{\psi}(\eta, \mathbf{b}) \in H^{2}_m(\mathbb{C}),
\end{equation}
be nontrivial solutions of (\ref{fulleq-5}) proved in Lemma \ref{lem-per-arg}. Here $(\eta, \mathbf{b})$ belongs to some open set $\mathcal{U}$ of $(0, \mathbf{0})$. We seek a criterion that characterizes a parameter regime of spiral patterns. Indeed, let $\widetilde{p}\,'(\eta, \mathbf{b})(s)$ be the phase derivative of $\widetilde{\psi}(\eta, \mathbf{b})$ in the polar form. Then the parameter regime of spiral patterns is given by
\begin{equation}
\{  (\eta, \mathbf{b}) \in \mathcal{U} : \widetilde{p} \,'(\eta, \mathbf{b})(s) \mbox{ is not identically zero.} \}.
\end{equation}

\begin{lemma}[spiral criterion] \label{lem-par-sub}
Suppose that $(\widetilde{\Omega}(\eta, \mathbf{b}), \widetilde{\psi}(\eta,\mathbf{b})) \in \mathbb{R} \times H_m^2(\mathbb{C})$ is a solution of \emph{(\ref{fulleq-5})}. If in the polar form of $\widetilde{\psi}(\eta,\mathbf{b})$ the phase derivative $\widetilde{p}\,'(\eta, \mathbf{b})(s)$ is identically zero, then
\begin{equation} \label{frorel}
\widetilde{\Omega}(\eta, \mathbf{b}) - \eta \, f_{\mathcal{R}}(0, \mathbf{b})  + f_{\mathcal{I}}(0,  \mathbf{b}) = 0.
\end{equation}
Consequently, $\{  (\eta, \mathbf{b}) \in \mathcal{U} : (\ref{frorel}) \emph{\mbox{ is violated.}}\}$ is a parameter subregime of spiral patterns.
\end{lemma}
\begin{proof}
Since $\widetilde{p}\,'(\eta, \mathbf{b})(s)$ is identically zero, $s = 0$ is an isolated zero of $\widetilde{A}(s)= |\widetilde{\psi}(\eta,\beta)(s)|$ (see Lemma \ref{asym} (i)), and $a(s) > 0$ for $s \in (0,s_*)$, by (\ref{trick}) there is a constant $\delta > 0$ such that 
\begin{equation} \label{threepara}
\widetilde{\Omega}(\eta,\mathbf{b}) - \eta \, f_{\mathcal{R}}(|\widetilde{A}(s)|^2, \mathbf{b}) + f_{\mathcal{I}}(|\widetilde{A}(s)|^2, \mathbf{b}) = 0 \quad \mbox{for all   } s \in (0, \delta).
\end{equation}
Since $\lim_{s \searrow 0}\widetilde{A}(s) = 0$, continuity implies (\ref{frorel}) as we take $s \searrow 0$ in (\ref{threepara}).
\end{proof}

The following lemma is equivalent to Theorem \ref{thm-gen} (i).

\begin{lemma} \label{final-1} 
Let $\lambda_0 > \lambda_{0}^{m}$ be fixed and $(\widetilde{\Omega}(\eta, \mathbf{b}), \widetilde{\psi}(\eta,\mathbf{b}))$ be nontrivial solution pairs of \emph{(\ref{fulleq-5})}. Then there exists an $\epsilon>0$ such that the following statements hold as $(\eta, \mathbf{b}) \in B_{\epsilon}^{d+1}(\emph{\textbf{0}})$:
\begin{itemize}
\item[(i)] There exists a smooth function $\widetilde{\eta} = \widetilde{\eta}(\mathbf{b})$, $\widetilde{\eta}: B^d_{\epsilon}(\emph{\textbf{0}}) \rightarrow \mathbb{R}$ with $\widetilde{\eta}(\emph{\textbf{0}}) = 0$ such that $\widetilde{\Omega}(\eta, \mathbf{b}) = 0$ if and only if $\eta = \widetilde{\eta}(\mathbf{b})$. 
\item[(ii)] There exists a smooth function $\kappa_1 = \kappa_1(\eta)$, $\kappa_1:(-\epsilon, \epsilon) \rightarrow [0, \infty)$ such that $\kappa_1(\eta) > 0$ if $\eta \neq 0$, and $( \widetilde{\Omega}(\eta, \mathbf{b}), \widetilde{\psi}(\eta, \mathbf{b}))$ exhibits a rotating spiral for all $\eta \neq 0$ and $0 \le |\mathbf{b}| \le \kappa_1(\eta)$.
\end{itemize}
\end{lemma}
\begin{proof}
For (i), the main idea is to use the frequency-parameter relation (\ref{fre-par-rel}),
\begin{equation} \label{fre-par-rel-5}
\mathcal{J} (\widetilde{\Omega}(\eta, \mathbf{b}), \widetilde{\psi}(\eta, \mathbf{b}), \, \eta, \mathbf{b}) := \int_\mathcal{M} \left( \widetilde{\Omega} - \eta  \, f_\mathcal{R}(|\widetilde{\psi}|^2, \mathbf{b}) + f_\mathcal{I}(|\widetilde{\psi}|^2, \mathbf{b}) \right) |\widetilde{\psi}|^2 = 0,
\end{equation}
evaluated at $(\eta, \mathbf{b}) \in B_{\epsilon}^{d+1}(\textbf{0})$. Using $\widetilde{\Omega}(0,\textbf{0}) = 0$, $\widetilde{\psi}(0,\textbf{0}) = \psi_0$, and the assumption $f_\mathcal{I}(y;\textbf{0}) = 0$ for all $y \ge 0$, we differentiate $\mathcal{J}$ with respect to $\eta$, evaluate at $(\eta, \mathbf{b}) = (0, \textbf{0})$, and then obtain
\begin{equation} \label{keyrel}
\partial_\eta \widetilde{\Omega}(0,\textbf{0}) = \frac{\int_\mathcal{M} f_\mathcal{R}(|\psi_0|^2,\textbf{0}) \, |\psi_0|^2}{\int_\mathcal{M} |\psi_0|^2}.
\end{equation}

Since $\psi_0$ satisfies the $C^0$-bound (\ref{priori}), (A1--A2) and (\ref{keyrel}) imply
\begin{equation} \label{open1}
0 < \partial_\eta \widetilde{\Omega}(0,\textbf{0}) < 1.
\end{equation}
 
Since $\widetilde{\Omega}(0, \textbf{0}) = 0$ and $0 < \partial_\eta \widetilde{\Omega}(0,\textbf{0})$, the implicit function theorem yields the existence of $\widetilde{\eta}$. This proves (i)

To prove (ii), by the spiral criterion (\ref{frorel}), we note that $( \widetilde{\Omega}(\eta, \mathbf{b}), \widetilde{\psi}(\eta,\mathbf{b}))$ exhibits a spiral if 
\begin{equation} \label{criteria-4}
\widetilde{\Omega}(\eta, \mathbf{b}) - \eta \, f_\mathcal{R}(0,\mathbf{b}) + f_\mathcal{I}(0,\mathbf{b}) \neq 0.
\end{equation}
Since $f_\mathcal{R}(0;\textbf{0}) = 1$, we have $\partial_\eta \widetilde{\Omega}(0,\textbf{0}) < f_\mathcal{R}(0;\textbf{0})$ by (\ref{open1}). Since $\widetilde{\Omega}(0,\textbf{0}) = 0$, smoothness of $\widetilde{\Omega}$ implies $\widetilde{\Omega}(\eta, \textbf{0}) \neq \eta \, f_\mathcal{R}(0;\textbf{0})$ for all $\eta \neq 0$ sufficiently near zero. Since $f_\mathcal{I}(0, \mathbf{0}) = 0$, smoothness of $f_\mathcal{R}$ and $f_\mathcal{I}$ yields a desired smooth function $\kappa_1 = \kappa_1(\eta)$. 
\end{proof}

Notice that Lemma \ref{final-1} (ii) does not assure the existence of spiral patterns for $\eta$ being almost zero, but $|\mathbf{b}|$ being large. Such a result of existence requires more assumptions on $f_\mathcal{I}$. For instance, recall the two assumptions in Theorem \ref{thm-gen} (ii),
\begin{equation} \label{hyp3}
\partial_{\beta} f_\mathcal{I}(y,0) \neq 0 \quad \mbox{for all   } y \in (0, C),
\end{equation}
\begin{equation} \label{hyp4}
f_\mathcal{I}(0,\beta) = 0 \quad \mbox{for all   } \beta \in \mathbb{R}.
\end{equation}

\begin{remark}
The assumptions (\ref{hyp3}) and (\ref{hyp4}) include the nonlinearity considered in \cite{KoHo81, Ts10}, which is given by $f_\mathcal{I}(y, \beta) = \beta \, \omega(y)$, where $\omega(0) = 0$ and $\omega(y) \neq 0$ for all $y \in (0, C)$.
\end{remark}

The following lemma implies Theorem \ref{thm-gen} (ii).

\begin{lemma}  \label{lem-d=1}
Suppose $d =1$ and under the same assumptions of Lemma \emph{\ref{final-1}}. Then the following statements hold:
\begin{itemize}
\item[(i)] Assume \emph{(\ref{hyp3})}. Then $\widetilde{\eta} = \widetilde{\eta}(\beta)$ proved in Lemma \ref{final-1} (ii) is invertible.
\item[(ii)] Assume \emph{(\ref{hyp3}--\ref{hyp4})}. Then there exists a smooth function $\kappa_2 = \kappa_2(\beta)$, $\kappa_2: (-\epsilon, \epsilon) \rightarrow [0, \infty)$ such that $\kappa_2(\beta) > 0$ if $\beta \neq 0$, and $( \widetilde{\Omega}(\eta, \beta), \widetilde{\psi}(\eta, \beta))$ exhibits a rotating spiral for all $\beta \neq 0$ and $0 \le |\eta| \le \kappa_2(\beta)$. 
\end{itemize}
\end{lemma}
\begin{proof}
For (i), we differentiate $\mathcal{J}$ in (\ref{fre-par-rel-5}) with respect to $\beta$ and evaluate at $(\eta, \beta) = (0, 0)$. Then (\ref{hyp3}) implies 
\begin{equation} \label{criteria-3}
\partial_{\beta}\widetilde{\Omega}(0,0) =  \frac{-\int_\mathcal{M} \partial_{\beta} f_\mathcal{I}(|\psi_0|^2,0) \, |\psi_0|^2}{\int_\mathcal{M} |\psi_0|^2} \neq 0.
\end{equation}
Thus $\widetilde{\eta}$ is invertible. 

For (ii), by (\ref{hyp4}) the spiral criterion (\ref{criteria-4}) now reads
\begin{equation} \label{criteria-6}
\widetilde{\Omega}(\eta, \beta) - \eta \, f_\mathcal{R}(0;\beta) \neq 0.
\end{equation}
Due to $\widetilde{\Omega}(0,0) = 0$ and (\ref{criteria-3}), we have $\widetilde{\Omega}(0, \beta) \neq 0$ for all $\beta \neq 0$ sufficiently near zero. Then smoothness of $f_\mathcal{R}$ yields a desired smooth function $\kappa_2 = \kappa_2(\beta)$. 
\end{proof}

For general complex GLe, the local information of $(\eta, \mathbf{b}) = (0, \textbf{0})$ is not sufficient for the spiral criterion (\ref{criteria-4}) to determine the type of pattern for all parameters $(\eta, \beta)$ in $B_{\epsilon}^{d+1}(\textbf{0})$. Interestingly, the cubic supercritical GLe gives a complete classification. 

The following lemma is equivalent to Theorem \ref{thm-cgle}.

\begin{lemma} \label{final-2}
Under the same assumptions of Lemma \emph{\ref{final-1}}, suppose $\mathbf{b} = \beta \in \mathbb{R}$ and $f(y,\beta) = 1- y - i \,\beta \,  y$. Then the following statements hold as $(\eta, \beta) \in B_{\epsilon}^{2}(\emph{\textbf{0}})$:
\begin{itemize}
\item[(i)] $\widetilde{p}\,'(\eta,\beta)(s)$ is identically zero if and only if $\eta = \beta$.
\item[(ii)] The smooth function $\widetilde{\eta} = \widetilde{\eta}(\beta)$ proved in Lemma \emph{\ref{final-1}} (ii) is strictly decreasing. 
\end{itemize}
Consequently, $(\widetilde{\Omega}(\eta, \beta), \widetilde{\psi}(\eta, \beta))$ exhibits a rotating spiral if $\eta \neq \beta$ and $\eta \neq \widetilde{\eta}(\beta)$; it exhibits a frozen spiral if $\eta \neq \beta$ and $\eta = \widetilde{\eta}(\beta)$. 
\end{lemma}
\begin{proof}
For (i), suppose that $\widetilde{p}\,'(\eta, \beta)(s)$ is identically zero. The spiral criterion (\ref{criteria-4}) and the specific form of $f_{\mathcal{I}}$ imply $\widetilde{\Omega}(\eta, \beta) = \eta$. Therefore, the frequency-parameter relation (\ref{fre-par-rel-5}) reads $\int_\mathcal{M} (  \eta - \beta) \, |\widetilde{\psi}|^4 = 0$, and so $\eta = \beta$. Conversely, suppose $\eta = \beta$. Then (\ref{fre-par-rel-5}) reads $\int_\mathcal{M} (  \widetilde{\Omega}(\eta,\eta) - \eta ) \, |\widetilde{\psi}|^2 = 0$, and so $\widetilde{\Omega}(\eta, \eta) = \eta$. Hence the spiral Ansatz is of the form $\Psi(t, s, \varphi) = e^{-i \eta t} \, \widetilde{\psi}(s,\varphi)$. Since $\Psi$ solves
\begin{equation}
\partial_t \Psi = \frac{1}{\lambda_0}(1 + i\, \eta) \Delta_{\mathcal{M}} \, \Psi + ( 1 - |\Psi|^2 - i\, \eta \, |\Psi|^2) \, \Psi, 
\end{equation}
it is easy to verify that $\widetilde{\psi}$ solves the following real GLe:
\begin{equation}
\Delta_{m} \widetilde{\psi} + \lambda_0\, ( 1 - |\widetilde{\psi}|^2 ) \, \widetilde{\psi} = 0. 
\end{equation}
Therefore, $\widetilde{p}\,'(\eta,\eta)(s)$ is identically zero due to the decoupling effect; see Lemma \ref{lem-dec-eff}. This proves (i).

For (ii), since $f_\mathcal{I}(y,\beta) = - \beta \,y$ fulfills the assumption (\ref{hyp3}), $\widetilde{\eta} = \widetilde{\eta}(\beta)$ is invertible by Lemma \ref{lem-d=1} (i). It remains to show $D_\beta \widetilde{\eta}(0) < 0$. Differentiating the relation $\widetilde{\Omega}(\widetilde{\eta}(\beta), \beta) = 0$ with respect to $\beta$ and evaluating at $\beta = 0$ yield
\begin{equation}
D_\beta \widetilde{\eta}(0) = - \frac{\partial_{\beta}\widetilde{\Omega}(0,0)}{\partial_{\eta}\widetilde{\Omega}(0,0)}.
\end{equation}
By (\ref{open1}) we have $\partial_{\eta}\widetilde{\Omega}(0,0) > 0$. Since $\partial_{\beta} f_\mathcal{I}(y;0) = -y$, by (\ref{criteria-3}) we see $\partial_{\beta}\widetilde{\Omega}(0,0) > 0$. Hence $D_\beta \widetilde{\eta}(0) < 0$. The proof is complete.
\end{proof}

\section{Discussion and outlook} \label{sec;6}

The Ginzburg-Landau spiral waves we have obtained possibly possess two features: slowly rotating and slightly twisting, due to perturbation arguments. Such spirals differ from those observed in the diffusive Belousov-Zhabotinsky equation; see \cite{Beetal97}. These two features are rooted in the global $S^1$-equivariance (\ref{glo-gau}), and the assumption (\ref{fci}) by which the Ginzburg-Landau equation is weakly coupled for  small parameters $0 < |\eta|,\, |\mathbf{b}| \ll 1$. To seek fast rotating or greatly twisting spiral waves, we shall either consider sufficiently large parameters (see the heuristic attempts in \cite{Ha82}) or replace (\ref{fci}).

We indicate two directions of research regarding our result of existence.

First, we can seek Ginzburg-Landau scroll waves, say on the three-dimensional unit ball. The major difference is that the resulting elliptic equation is not equivalent to a second-order ODE, and thus the principal bifurcation curve may not possess nodal structure (see Lemma \ref{lem-nod-4}). Therefore, new ideas are needed to extend the bifurcation curve globally.

Second, we shall rigorously analyze the stability of spiral waves. Indeed, formal asymptotic expansions and numerical evidences suggest that the one-armed spiral waves be stable, while multi-armed ones be unstable; see \cite{Ha82, Ts10}. Furthermore, we can study whether the stability changes when noninvasive spatio-temporal feedback control is introduced, by the method established in \cite{Sch16}.



\section*{Acknowledgments}

This work was supported by the National Center for Theoretical Sciences through grant 107-2119-M-002-016. The author is grateful to Bernold Fiedler for his PhD supervision, and all members of the Nonlinear Dynamics group at the Free University of Berlin.

\bibliographystyle{siamplain}
\bibliography{references}
\end{document}

%% file: ex_shared.tex
\usepackage{lipsum}
\usepackage{amsfonts}
\usepackage{graphicx}
\usepackage{epstopdf}
\usepackage{algorithmic}
\ifpdf
  \DeclareGraphicsExtensions{.eps,.pdf,.png,.jpg}
\else
  \DeclareGraphicsExtensions{.eps}
\fi

\newsiamremark{remark}{Remark}
\newsiamremark{hypothesis}{Hypothesis}
\crefname{hypothesis}{Hypothesis}{Hypotheses}
\newsiamthm{claim}{Claim}

\headers{Ginzburg-Landau Spiral Waves in Circular and Spherical Geometries}{Jia-Yuan Dai}

\title{Ginzburg-Landau Spiral Waves in Circular and Spherical Geometries\thanks{Submitted to the editors DATE.
\funding{This work was funded by the National Center for Theoretical Sciences through grant 107-2119-M-002-016.}}}

\author{Jia-Yuan Dai\thanks{National Center for Theoretical Sciences, National Taiwan University, Taiwan (\email{jydai@ncts.ntu.edu.tw}).}
}

\usepackage{amsopn}
